\documentclass[12pt,reqno]{amsart}

\usepackage{graphicx}
\usepackage[all]{xy}

\theoremstyle{plain}
  \newtheorem{theorem}{Theorem}
  \newtheorem{corollary}{Corollary}
  \newtheorem{lemma}{Lemma}
  \newtheorem{proposition}{Proposition}
\theoremstyle{definition}
  \newtheorem{definition}{Definition}
  \newtheorem*{question*}{Question}
  
  \newtheorem*{todo*}{TODO}
\theoremstyle{remark}
  \newtheorem{remark}{Remark}

\newcommand\itemref[1]{(\ref{#1})}

\newcommand\ie{i.e.\ }
\newcommand\cf{cf.\ }
\newcommand\resp{resp.\ }

\newcommand\an{\text{\rm an}}
\newcommand\odd{\text{\rm odd}}
\newcommand\even{\text{\rm even}}
\newcommand\op{\text{\rm op}}
\newcommand\geom{\text{\rm geom}}
\newcommand\met{\text{\rm met}}
\newcommand\comb{\text{\rm comb}}

\newcommand\Hodge{\text{\rm Hodge}}

\newcommand\cS{\mathcal S}
\newcommand\cR{\mathcal R}

\newcommand\Eul{\mathfrak{Eul}}

\newcommand\Z{\mathbb Z}
\newcommand\N{\mathbb N}
\newcommand\R{\mathbb R}
\newcommand\K{\mathbb K}
\newcommand\X{\mathbb X}

\newcommand\M{\mathcal M}
\newcommand\T{\mathcal T}

\newcommand\e{\mathfrak e}

\DeclareMathOperator{\inc}{inc}
\DeclareMathOperator{\pr}{pr}

\DeclareMathOperator{\pt}{pt}
\DeclareMathOperator{\tr}{tr}

\DeclareMathOperator{\Or}{\mathcal O}
\DeclareMathOperator{\rank}{rank}
\DeclareMathOperator{\img}{img}
\DeclareMathOperator{\IND}{IND}
\DeclareMathOperator{\ind}{ind}

\DeclareMathOperator{\id}{id}

\DeclareMathOperator{\grad}{grad}
\DeclareMathOperator{\cs}{cs}
\DeclareMathOperator{\Det}{Det}
\DeclareMathOperator{\Vol}{vol}
\DeclareMathOperator{\Int}{Int}

\begin{document}

\title[The geometric complex of a Morse--Bott--Smale pair \dots]
      {The geometric complex of a Morse--Bott--Smale pair
       and an extension of a theorem by Bismut--Zhang}

\author{Dan Burghelea}

\address{Dan Burghelea,
         Dept. of Mathematics,
         The Ohio State University,
         231 West Avenue, Columbus, OH 43210, USA.}

\email{burghele@mps.ohio-state.edu}

\author{Stefan Haller}

\address{Stefan Haller,
         Department of Mathematics, University of Vienna,
         Nordbergstra{\ss}e 15, A-1090, Vienna, Austria.}

\email{stefan.haller@univie.ac.at}

\thanks{Part of this work was done while both authors enjoyed the 
        hospitality of the Max Planck Institute for Mathematics.
        The second author is supported by the \emph{Fonds zur 
        F\"orderung der wis\-sen\-schaft\-lichen Forschung} (Austrian 
        Science Fund), project number {\tt P14195-MAT}}



\date{\today}

\begin{abstract}
  We define the geometric complex associated to a Morse--Bott--Smale
  vector field, \cf \cite{AB95}, and its associated spectral sequence.
  We prove an extension of the Bismut--Zhang theorem
  to Morse--Bott--Smale functions, see Theorem~\ref{T:intro} below. 
  The proof is based on the Bismut--Zhang theorem for Morse--Smale functions, see
  \cite{BZ92}.
\end{abstract}

\maketitle

\setcounter{tocdepth}{1}
\tableofcontents

  \section{Introduction}
  \label{S:intro}

Let $M$ be a closed connected manifold, $g$ a Riemannian metric on $M$,
$E$ a flat vector bundle over $M$ and $\mu$ a fiber metric on $E$.
This permits to define the formal adjoint of the deRham differential
on $\Omega^*(M;E)$ and the Laplacian $\Delta^q$ acting on $\Omega^q(M;E)$.
In this situation the Ray--Singer torsion \cite{RS71} is defined as:
$$
\log T_\an^M(E,g,\mu)
:=\frac12\sum_q(-)^{q+1}q\log\det{}'\Delta^q
$$
where $\det{}'\Delta^q$ denotes the zeta regularized determinant of 
$\Delta^q$ with zero modes discarded.
The main motivation of this paper is the calculation of the Ray--Singer torsion.
This will be done with the help of the geometric complex associated  to a 
Morse--Bott--Smale vector field described below, which is an object
of independent interest by itself.

Suppose $X$ is a Morse--Bott--Smale vector field on $M$, \ie 
$X=-\grad_{\tilde g}(f)$ where $(f,\tilde g)$ is a Morse--Bott--Smale pair, \cf
Definition~\ref{D:MBS} in section~\ref{SS:traj}. Let $\Sigma$
denote the critical manifold of $X$ and $N^-\to\Sigma$ the negative normal 
bundle. For every connected component $S\subseteq\Sigma$ we have a Morse 
index $\ind(S)\in\N_0:=\N\cup\{0\}$ which coincides with the rank of $N^-$ over $S$. 
Every connected component $S\subseteq\Sigma$ inherits a Riemannian metric
$g_S$ from $g$. Let $\Or_{N^-}$ denote the orientation bundle of $N^-$.
The flat vector bundle $E_S:=(E|_\Sigma\otimes\Or_{N^-})|_S$ 
over $S$ inherits a fiber metric $\mu_S$ from the fiber metric
$\mu|_S$ on $E|_S$ and the canonic parallel fiber metric on
$\Or_{N^-}$.\footnote{The frame bundle of $N^-$ induces a principal
$\Z_2$--covering of $\Sigma$. We think of the orientation bundle $\Or_{N^-}$
as the flat line bundle associated to this principal covering with respect to the
standard non-trivial $\Z_2$--action on $\R$. It thus inherits a canonic
parallel fiber metric from the standard Euclidean inner product on $\R$.} 
So we have Ray--Singer torsions 
$\log T_\an^S(E_S,g_S,\mu_S)$ for every connected component 
$S\subseteq\Sigma$.

The Morse--Bott--Smale vector field $X$ gives rise to a geometric complex
$C^*(X;E)$ whose underlying graded vector space is
$$
C^*(X;E):= 
\bigoplus_{S\subseteq\Sigma}\Omega^{*-\ind(S)}(S;E_S)
$$  
where the sum is over all connected components $S\subseteq\Sigma$.
Fiber integration over the unstable manifolds provides a 
homomorphism of complexes $\Int:\Omega^*(M;E)\to C^*(X;E)$, see
section~\ref{S:comp} below. 
The Morse index provides a decreasing finite filtration 
$$
C^*_p(X;E):=\bigoplus_{\ind(S)\geq p}\Omega^{*-\ind(S)}(S;E_S)
$$
on the geometric complex. This gives rise to a spectral sequence 
$(E_kC^*(X;E),\delta_k)$, see for instance \cite{BT82}.

\begin{theorem}\label{T:thom}
\
\begin{enumerate}
\item\label{T:thom:i}
The integration map induces an isomorphism in cohomology
$\Int:H^*(M;E)=HC^*(X;E)$.
\item\label{T:thom:ii}
The spectral sequence $(E_kC^*(X;E),\delta_k)$ converges to $HC^*(X;E)$.
Moreover $E_1C^*(X;E)=\bigoplus_{S\subseteq\Sigma}H^{*-\ind(S)}(S;E_S)$.
Particularly $E_kC^*(X;E)$ is finite dimensional for all $k\geq1$
and we have $E_kC^*(X;E)=E_\infty C^*(X;E)$ for sufficiently large $k$.
\end{enumerate}
\end{theorem}

The above theorem was first stated in \cite{AB95} but known to experts. 
Unfortunately, despite its usefulness and simplicity there is no complete
(satisfactory) treatment  of the geometric complex in the existing 
literature. Since we need the geometric complex for our main result  
we will complete the discussion of \cite{AB95}
and provide a complete treatment of the geometric complex and its basic
properties.

\begin{remark}
If $X$ is a Morse--Smale vector field the geometric complex is finite 
dimensional, $(E_1C^*(X;E),\delta_1)=(C^*(X;E),\delta)$ is the familiar Morse
complex, and $\delta_k=0$ for $k\geq2$. In general $(E_1C^*(X;E),\delta_1)$ 
has a very simple form, see section~\ref{SS:thom}.
\end{remark}

\begin{remark}
If $\pi:M\to N$ is a smooth bundle with $M$ and $N$ closed manifolds, 
and $X_N$ is a Morse--Smale vector field on $N$ then one can construct a 
Morse--Bott--Smale vector field on $X_M$ on $M$ with
$d_x\pi(X_M(x))=X_N(\pi(x))$ for all $x\in M$. 
The spectral sequence associated with the geometric complex of 
$X_M$ identifies with the Leray--Serre spectral sequence of the bundle 
$\pi$ provided by the cell structure on $N$ induced by $X_N$, \cf 
Propositions~\ref{P:MBS_on_fb} and \ref{P:spec_for_fb} in
section~\ref{SS:bundles}.
\end{remark}

The Riemannian metric $g$ and the fiber metric $\mu$ induce inductively  
scalar products on every $E_k$--term of the spectral sequence. Indeed,
the restriction of $g$ and $\mu$ to $\Sigma$ provide a scalar product
on the deRham complex $\Omega^*(\Sigma;E_\Sigma)$. Its restriction to harmonic
forms provides the scalar product on $E_1C^*(X;E)$, \cf
Theorem~\ref{T:thom}\itemref{T:thom:ii}. Now suppose inductively that
$E_{k-1}C^*(X;E)$ is equipped with a scalar product. Then
$E_kC^*(X;E)=H(E_{k-1}C^*(X;E),\delta_{k-1})$ inherits a scalar product
as a subquotient of the finite dimensional Euclidean vector space.
Equivalently, one can describe this scalar product via
finite dimensional Hodge theory on the complex $(E_{k-1}C^*(X;E),\delta_k)$.

The scalar products on $E_kC^*(X;E)$ permit to define the adjoint of the 
differential $\delta_k$ on $E_kC^*(X;E)$ and a Laplacian $\Delta_k^q$ acting 
on the finite dimensional $E_kC^q(X;E)$, $k\geq1$. Define the combinatorial
torsion by:
$$
\log T_\comb^M(E,X,g,\mu)
:=\sum_{k\geq1}\frac12\sum_q(-)^{q+1}q\log\det{}'\Delta^q_k
$$

The deRham cohomology inherits a scalar
product via Hodge decomposition which provides a scalar product on the one 
dimensional real vector space $\det H^*(M;E)$ referred to as the 
\emph{Hodge scalar product}.\footnote{Here we use the notation
$\det V=\Lambda^{\dim(V^\even)}V^\even\otimes\Lambda^{\dim(V^\odd)}(V^\odd)^*$
for a finite dimensional $\Z_2$--graded vector space. See appendix~\ref{A:homalg} 
for more details on determinant lines.} 
The canonical identification $\det HC^*(X;E)=\det E_{\infty}C^*(X;E)$ and
the scalar product on $E_\infty C^*(X;E)=E_kC^*(X;E)$ for $k$ large,
define a scalar product on $\det HC^*(X;E)$ referred to as the 
\emph{geometric scalar product}.
The isomorphism in cohomology induced by integration provides an isomorphism 
$\det\Int:\det H^*(M;E)\to\det HC^*(X;E)$. Define\footnote{Recall that 
for an isomorphism $\alpha:V\to W$ between Euclidean vector
spaces one defines $\Vol\alpha:=\sqrt{\det(\alpha^*\alpha)}$.}
$$
\log T_\met^M(E,X,g,\mu):=\log\Vol(\det\Int).
$$

The fiber metric $\mu$ on $E$ provides a closed one form
$\theta\in\Omega^1(M;\R)$ which measures to what extend the induced
length one section on $\det E\otimes\Or_E$ is not parallel, see 
section~\ref{SS:an_tor}. Define
$$
\cR(\theta,X,g):=\int_{M\setminus\Sigma}\theta\wedge X^*\Psi_g
$$
where $\Psi_g\in\Omega^{\dim M-1}(TM\setminus M;\Or_M)$ denotes the 
global angular form, \cf \cite{BT82}, also called Mathai--Quillen form 
in \cite{BZ92}. The integral on the right hand side need not be 
convergent but can be regularized, see section~\ref{SS:euler}.

The main theorem of this paper is the following extension
of a result by Bismut--Zhang \cite{BZ92},
Cheeger \cite{Ch77}, \cite{Ch79} and M\"uller \cite{Mu78}.

\begin{theorem}\label{T:intro}
Let $(M,g)$ be a closed Riemannian manifold, $E$ a flat vector bundle with
fiber metric $\mu$ over $M$, and let $X$ be a Morse--Bott--Smale vector 
field on $M$, see Definition~\ref{D:MBS} below. Then
\begin{multline*}
\log T_\an^M(E,g,\mu)
=\sum_{S\subseteq\Sigma}(-)^{\ind(S)}\log T_\an^S(E_S,g_S,\mu_S)
\\
+\log T_\comb^M(E,X,g,\mu)
+\log T_\met^M(E,X,g,\mu)
+\cR(\theta,-X,g)
\end{multline*}
where the sum is over all connected components $S\subseteq\Sigma$.
\end{theorem}

One can regard $\Sigma$ as a graded manifold
where the degree of the component $S$ is $\ind(S)$ and write 
$\log T_\an^\Sigma:=\sum_{S\subseteq\Sigma}(-)^{\ind(S)}\log T_\an^S$.
Then one can view the theorem above as a localization theorem stating that
$$
\log T_\an^M=\log T_\an^\Sigma+\log T_\comb+\log T_\met+\cR(\theta,-X,g)
$$
with the term $\cR(\theta,-X,g)$ computable in terms of differential 
geometry; the term $\log T_\comb$ computable in terms of 
finite dimensional linear algebra; and, in case $H^*(M;E)=0$, vanishing term 
$\log T_\met$.

Theorem~\ref{T:intro} can be extended to the case $E$ is a flat bundle of
Hilbertian $\mathcal A$--modules of finite type where $\mathcal A$ is a 
finite von Neuman algebra using the extension of the 
Bismut--Zhang theorem to this case as formulated in \cite{BFK01}. 
The proof and the results are exactly the same.

Theorem~\ref{T:intro} can be extended to $G$--manifolds where $G$ is a
finite group and the torsions replaced by the $G$--equivariant torsions
which are class functions on $G$ with positive real values.
This will use the corresponding extension of the Bismut--Zhang 
theorem \cf \cite{BZ94}. The particular case $G=\Z_2$ when applied to 
the double of a manifold with boundary \resp to a bordism will imply an 
extension of  Theorem~\ref{T:intro} to compact manifolds with boundary 
\resp to bordisms and, as application, will provide a very simple derivation 
of all known glueing formulae for analytic torsion \cf \cite{BFK99}. This is 
contained in \cite{BH04}.

Theorem~\ref{T:intro} permits to recover and generalize results of 
L\"uck--Schick--Thielmann \cf \cite{LST98} about the torsion of smooth 
bundles, \cf section~\ref{S:bundles}. In a forthcoming paper we will use the 
above results and reduce the calculation of $G$--equivariant torsion for 
$G$ a compact Lie group to finite dimensional linear algebra.

There is a more conceptual way to formulate Theorem~\ref{T:intro} and
understand the extend to which the result is a localization theorem for
torsion. To explain this consider $M$ a closed manifold, $E$ a flat vector
bundle over $M$, $x_0\in M$ a base point, $X$ a Morse--Bott--Smale vector
field which gives rise to the geometric complex $C^*(X;E)$, and
$\e\in\Eul_{x_0}(M)$ an Euler structure based at $x_0$, see
section~\ref{SS:euler}. The analytic torsion $\tau^{E,\e}_\an$ can be
regarded as an element defined up to multiplication by $\pm1$ in the one
dimensional vector space\footnote{Here we use the following convention.
For a line $L$, \ie a one dimensional vector space, we write 
$L^k:=L\otimes\cdots\otimes L$,
$L^{-k}=L^*\otimes\cdots\otimes L^*$ and we write $V^0:=\R$.
Note that there is a canonic isomorphism $L^{k_1+k_2}=L^{k_1}\otimes L^{k_2}$
for all $k_1$ and $k_2$.}
$$
\Det_{x_0}(M;E):=\det H^*(M;E)\otimes(\det E_{x_0})^{-\chi(M)}
$$
where $H^*(M;E)$ denotes the deRham cohomology with coefficients in $E$
and $\det H^*(M;E)$ its determinant line in the graded sense, see
appendix~\ref{A:homalg}. The geometric torsion $\tau_\geom^{E,\e,X}$ is an
element well defined up to multiplication by $\pm1$ in the one dimensional
vector space
$$
\Det_{x_0}(X;E):=\det HC^*(X;E)\otimes(\det E_{x_0})^{-\chi(M)}
$$
where $HC^*(X;E)$ denotes the cohomology of the geometric complex, see
section~\ref{S:geom_complex} for definition. The integration induces an 
isomorphism
$$
\det\Int:\Det_{x_0}(M;E)\to\Det_{x_0}(X;E)
$$
and Theorem~\ref{T:intro} is equivalent to the statement that
$\det\Int(\tau_\an^{E,\e})=\tau^{E,\e,X}_\geom$, \cf Theorem~\ref{T:main}
below.

The geometric torsion $\tau_\geom^{E,\e,X}$ is obtained from the analytic
torsions $\tau_\an^{E_S,\e_S}$ for a collection of Euler structures
$\e_S\in\Eul_{x_S}(S)$, one for every connected component $S\subseteq\Sigma$.
Precisely, the geometric complex $C^*(X;E)$ provides an isomorphism:
$$
I_1:\det HC^*(X;E)\to\bigotimes_{S\subseteq\Sigma}
\bigl(\det H^*(S;E_S)\bigr)^{(-)^{\ind(S)}}
$$
The Euler structure $\e$ and the Euler structures $\e_S$ provide a one chain
$c\in C_1(M;\R)$ which one writes as
$c=c_0+\sum_{S\subseteq\Sigma}(-)^{\ind(S)}\chi(S)c_S$
where $c_S$ is a path from $x_0$ to $x_S$ and $c_0\in C_1(M;\R)$ which then
satisfies $\partial c_0=0$. Parallel transport along the $c_S$ and proper scaling 
determined by $c_0$ provides an isomorphism well defined up to sign:
$$
I_2:(\det E_{x_0})^{-\chi(M)}\to\bigotimes_{S\subseteq\Sigma}
\bigl((\det E_S)_{x_S}^{-\chi(S)}\bigr)^{(-)^{\ind(S)}}
$$
The geometric torsion is then defined by:
$$
\tau_\geom^{E,\e,X}:=
(I_1\otimes I_2)^{-1}\Bigl(
\bigotimes_{S\subseteq\Sigma}\bigl(\tau_\an^{E_S,\e_S}\bigr)^{(-)^{\ind(S)}}
\Bigr)
$$
This explains the localization aspect of the result.

The paper contains a number of intermediate results of independent interest
such as: the canonic compactification of the unstable sets and the space of
trajectories of a Morse--Bott--Smale vector field, see
Propositions~\ref{P:comp_T} and \ref{P:comp_W} and Theorems~\ref{T:comp_T}
and \ref{T:comp_W}; the invariant $\cR$ and the discussion of Euler
structures.

  \section{Analysis of the ODE}
  \label{S:comp}

\subsection{Morse--Bott functions}
\label{SS:MB_func}

Let $M$ be a closed manifold, \ie $M$ is compact and without boundary. 
Let $f:M\to\R$ be a smooth function and let $\Sigma:=\{x\in M\mid df_x=0\}$
denote the \emph{critical set of $f$}. Suppose $x$ is a 
critical point of $f$ and $Y,Z\in T_xM$. Then
$H(Y,Z):=Y\cdot df(Z)=Z\cdot df(Y)$ for $df$ is closed
and vanishes at $x$. So we obtain a fiber wise symmetric bilinear form $H$
on the vector bundle $TM|_\Sigma$ referred to as the 
\emph{Hessian of $f$}.\footnote{This bilinear form can be extended 
(in a non-canonical way) to a smooth bilinear form on $TM$ as
follows. Pick any linear torsion free connection $\nabla$ on 
$TM$ and set $H(Y,Z):=(\nabla df)(Y,Z)$. This will be 
symmetric for $df$ is a closed one form and $\nabla$ was 
supposed to be torsion free. It will coincide with our
previous definition over $\Sigma$ since $df$ vanishes along
$\Sigma$.} If $\Sigma\subseteq M$ happens to be a closed 
submanifold the Hessian will provide a smooth fiber wise
bilinear form on the smooth vector bundle $TM|_\Sigma$.
Moreover we obviously have $T\Sigma\subseteq\ker H$.
Therefore the Hessian induces a fiber wise
bilinear form on the normal bundle $N:=TM|_\Sigma/T\Sigma$
of $\Sigma$ which we will denote by $H$ too.

\begin{definition}[Morse--Bott functions]\label{D:MB}
A function $f:M\to\R$ is called \emph{Morse--Bott} if its
critical set $\Sigma$ is a closed submanifold of $M$ and the 
Hessian considered as fiber wise bilinear form 
on the normal bundle $N$ of $\Sigma$ is non-degenerate.
In this case $\Sigma$ will also be referred to as the
\emph{critical manifold} or the \emph{rest manifold}.
The index, \ie the number of negative eigen values, 
of the Hessian provides a locally constant
$\ind:\Sigma\to\N_0$ called the \emph{Morse index of $f$}.
Finally, define the \emph{critical manifold of index $q$} by
$\Sigma_q:=\bigsqcup_{\ind(S)=q}S$ where the disjoint union 
is over all connected components $S\subseteq\Sigma$ of index $q$.
\end{definition}

The following well known result totally
clarifies the local behavior of $f$ near its critical
submanifold.

\begin{proposition}[Morse lemma]\label{P:Morse_lemma}
Let $f$ be a Morse--Bott function on $M$, let $p_N:N\to\Sigma$ denote the 
normal bundle of its critical manifold $\Sigma$ and let $H$ denote the 
Hessian. Define a smooth function $\frac12H^2:N\to\R$ on the total space of
$N$ by $\frac12H^2(Y):=\frac12H(Y,Y)$. Then there exists a tubular neighborhood
$\varphi:N\to M$ such that $\varphi^*df=d(\frac12H^2)$ holds
in a neighborhood of the zero section $\Sigma\subseteq N$.
Particularly $\varphi^*f=\frac12H^2+p_N^*(f|_\Sigma)$ with $f|_\Sigma$
locally constant.
A tubular neighborhood like this is called a \emph{Morse chart} for $f$.
\end{proposition}

\subsection{Compatible Riemannian metrics}

Suppose $f$ is a Morse--Bott function on a close manifold $M$.
It is easy to see that there are $H$--orthogonal subbundles 
$N^+\oplus N^-=N$ such that 
$g_\pm:=\pm H|_{N^\pm}>0$. Certainly $H=g_+\ominus g_-$
with respect to this splitting. The bundle $N^+$
is called the \emph{positive normal bundle} and the bundle
$N^-$ is called the \emph{negative normal bundle} of $f$.
Note however, that these bundles are well defined only
up to isomorphisms of vector bundles.
For the rank of $N^-$ we have $\rank(N^-)=\ind:\Sigma\to\N_0$.

Suppose we have a Riemannian metric $g_\Sigma$ on
$\Sigma$ and linear connections $\nabla^\pm$ on $N^\pm$
such that $\nabla^\pm g_\pm=0$. Define a fiber metric
$g_N:=g_+\oplus g_-$ and a linear connection 
$\nabla^N:=\nabla^+\oplus\nabla^-$ on $N$. Then
$\nabla^Ng_N=0=\nabla^NH$. Define a Riemannian metric 
$\tilde g$ on the total space of $N$ using $g_N$, $g_\Sigma$
and $\nabla^N$ in the most obvious way. Clearly the restriction
of $\tilde g$ to $\Sigma$ coincides with $g_\Sigma$. 
On the normal bundle of $\Sigma\subseteq N$, which canonically
identifies to $N\to\Sigma$, the Riemannian metric $\tilde g$
induces a fiber metric and a linear connection which coincide
with $g_N$ and $\nabla^N$. The splitting of $N$ into positive
and negative $H$--eigen spaces using $g_N$ coincides with
$N^+\oplus N^-=N$. The second fundamental form of
$\Sigma\subseteq N$ vanishes. The endomorphism 
$-(\id_{N^+}\ominus\id_{N^-})$ of $N$ can be considered as a 
vertical vector field on the total space of $N$ which 
coincides with $X:=-\grad_{\tilde g}(\frac12H^2)$. The unstable
set of $\Sigma$, \ie the points in $N$ whose $X$ trajectory
departs at $\Sigma$, coincides with $N^-\subseteq N$ hence is 
a smooth submanifold of $N$. Similarly the stable set of 
$\Sigma$ coincides with the smooth submanifold 
$N^+\subseteq N$.

\begin{definition}[Compatible Riemannian metrics]\label{D:comp_g}
Let $f$ be a Morse--Bott function on $M$. A Riemannian metric 
$g$ on $M$ is called \emph{compatible with $f$} if there 
exists a Morse chart $\varphi:N\to M$ such that in addition to 
$\varphi^*df=d(\frac12H^2)$ we have $\varphi^*g=\tilde g$ 
on a neighborhood of $\Sigma\subseteq N$. Here $\tilde g$
is a Riemannian metric on the total space of $N$
constructed from an $H$--orthogonal splitting $N^+\oplus N^-=N$
for which $\pm H|_{N^\pm}>0$, a Riemannian metric on $\Sigma$ 
and linear connections on $N^\pm$ for which the restrictions 
of $H$ are parallel as explained in the previous paragraph.
\end{definition}

Note that a Riemannian metric $g$ on $M$ induces a Riemannian metric $g_\Sigma$ on
$\Sigma$, a fiber metric $g_N$ on the normal bundle $N=(T\Sigma)^\perp$, an orthogonal
splitting $N^+\oplus N^-=N$ and a smooth map (the exponential map) 
$\exp:N\to M$ from the total space $N$ of the normal
bundle of $\Sigma$ to $M$. If $g$ is compatible then the parallel
transport along a path in $\Sigma$ leaves $N$ and $N^\pm$ invariant and 
therefore define the connection $\nabla^N$, and $\nabla^\pm$.
The Riemannian metric $g_\Sigma$, the fiber metric $g_N$ and the
connection $\nabla^N$ induce a Riemannian metric on $N$.  If $g$ 
is compatible then $\exp$ is an isometry when restricted to a small
neighborhood of the zero section of the normal bundle. This can also be
taken as an equivalent definition of compatible metric with respect to $f$.

Note that given an arbitrary metric $g$ and an arbitrary neighborhood $U$
of $\Sigma$ one can produce a Riemannian metric $g'$ which is compatible for 
$f$ and agrees with $g$ on $\Sigma$ and on $M\setminus U$.

\subsection{Stable and unstable manifold}
\label{SS:unstable_mf}

Let $f$ be a Morse--Bott function on a closed manifold $M$ and let 
$g$ be a compatible Riemannian metric on $M$, see Definition~\ref{D:comp_g}.
Define a vector field $X:=-\grad_g(f)$. For a connected
component $S\subseteq\Sigma$ define its \emph{stable} \resp 
\emph{unstable set} by
$$
W^\pm_S:=\bigl\{x\in M\bigm|\lim_{t\to\pm\infty}\phi_t(x)\in S\bigr\}
$$
where $\phi_t$ denotes the flow of the vector field $X$ at
time $t$. $W^\pm_S$ is a smooth boundary-less submanifold of $M$
and we let $i^S_\pm:W^\pm_S\to M$ denote the inclusions. Moreover
we get smooth fiber bundles $p^S_\pm:W^\pm_S\to S$ diffeomorphic
to $N^\pm_S\to S$.
Define $W^\pm:=\bigsqcup_{S\subseteq\Sigma}W^\pm_S$ where the 
disjoint union is over all connected components $S$ of $\Sigma$.
The $i_\pm^S$ give rise to a smooth immersion $i_\pm:W^\pm\to M$ 
and the $p_\pm^S$ give rise to a smooth fiber bundle
$p_\pm:W^\pm\to\Sigma$ diffeomorphic to $N^\pm\to\Sigma$.

Let $V_{W^\pm}\to W^\pm$ denote the 
vertical (vector) bundle of $p_\pm:W^\pm\to\Sigma$. For its rank we find
$\rank(V_{W^-})=p_-^*\ind:W^-\to\N_0$. Here we write 
$p_-^*\ind$ for the composition $\ind\circ p_-:W^-\to\Sigma\to\N_0$.
The Lie transport along the normalized $\frac X{|X|}$ provides
an isomorphism of vector bundles $V_{W^\pm}=p_\pm^*N^\pm$ and
thus a canonic isomorphism of orientation bundles:
\begin{equation}\label{E:A}
\Or_{V_{W^\pm}}=p_\pm^*\Or_{N^\pm}
\end{equation}
Moreover, if $E$ is a flat
vector bundle over $M$ then parallel transport along $X$
provides a canonic isomorphism of flat vector bundles:
\begin{equation}\label{E:B}
i_\pm^*E=p_\pm^*E|_\Sigma
\end{equation}

\subsection{The space of trajectories}
\label{SS:traj}

Let $f$ be a Morse--Bott function on a closed manifold $M$
and let $g$ be a compatible Riemannian metric on $M$, \cf
Definition~\ref{D:comp_g}.
For $x\in\Sigma$ let $i_-^x:W^-_x\to M$ denote the
restriction of $i_-:W^-\to M$ to the fiber $W^-_x$ of
$p_-:W^-\to\Sigma$ over $x$.

\begin{definition}[Morse--Bott--Smale vector fields]\label{D:MBS}
A pair $(f,g)$ consisting of a Morse--Bott function $f$ 
and a compatible Riemannian metric $g$ on a closed manifold 
$M$ is called \emph{Morse--Bott--Smale pair} if the mappings
$i_+:W^+\to M$ and $i_-^x:W^-_x\to M$ are transversal for
all $x\in\Sigma$.\footnote{Given a Morse--Bott function $f$ 
one can ask if there exists a Riemannian
metric $g$ so that the pair $(f,g)$ is a Morse--Bott--Smale pair. Such a
function can be called `good' Morse--Bott function. All Morse functions
are good Morse--Bott functions but, as expected, not all Morse--Bott
functions are good. To decide when a Morse--Bott function is good is an
interesting problem. Here are a few examples:

1) If $f:N\to\R$ is a Morse function and $\pi:M\to N$ is a smooth
bundle with compact fibers then $f\circ\pi$ is a good Morse--Bott function.

2) If $M$ is a compact smooth $G$--manifold with $G$ a compact Lie group,
any normal $G$--Morse function $f:M\to\R$ (\ie for any critical
point the negative isotropy representation is trivial) is a good 
Morse--Bott function \cf \cite{BH04}. It is known
that normal $G$--Morse functions are $C^0$--dense in the space of all smooth
$G$--invariant functions in $C^0$--topology.

3) If $f$ is a Morse--Bott function so is $-f$ but it is possible that $f$ is
good and $-f$ is not.}
A vector field $X$ is called Morse--Bott--Smale if there exists a 
Morse--Bott--Smale pair $(g,f)$ such that $X=-\grad_g(f)$.
\end{definition}

Suppose $X$ is a Morse--Bott--Smale vector field on $M$.
Let $S\neq S'$ be two different connected components of $\Sigma$.
Then $\M(S,S'):=W^-_S\cap W^+_{S'}$ is a smooth submanifold 
of $M$. Define $\M:=\bigsqcup_{S\neq S'\subseteq\Sigma}\M(S,S')$
where the disjoint union is over all pairs of different connected 
components $S$ and $S'$ of $\Sigma$. 
Let $i:\M\to M$ denote
the obvious smooth immersion, $p_+:\M\to\Sigma$ the obvious
smooth mapping and let $p_-:\M\to\Sigma$ denote the obvious
smooth fiber bundle.

Let $V_\M\to\M$ denote the vertical bundle of $p_-:\M\to\Sigma$.
We have a short exact sequence of vector bundles: 
$$
0\to V_\M\to V_{W^-}|_\M\to TM|_\M/TW^+|_\M\to0
$$
As in section~\ref{SS:unstable_mf}, Lie transport along $\frac X{|X|}$
provides canonic isomorphisms
$V_{W^-}|_\M=p_-^*N^-$ and $TM|_\M/TW^+|_\M=p_+^*N^-$. So we  
obtain a short exact sequence of vector bundles:
$$
0\to V_\M\to p_-^*N^-\to p_+^*N^-\to0
$$
Particularly $\rank(V_\M)=p_-^*\ind-p_+^*\ind:\M\to\N_0$.
Moreover, this short exact sequence provides us with a 
canonic isomorphism of orientation bundles:\footnote{Throughout the
paper we use the following convention. If $0\to E\to F\to G\to0$ is a 
short exact sequence of vector bundles we use an orientation of $G_x$
followed by an orientation of $E_x$ to yield the compatible orientation 
of $F_x$, $x$ in the common base space. We refer to this as the $G_x$
first convention. The notation in \eqref{E:AB} should indicate that we 
deviate from this convention and use the $V_\M$ first convention here.}
\begin{equation}\label{E:AB}
(-)^{\rank(V_\M)\cdot\rank(p_+^*N^-)}:
\Or_{V_\M}\otimes p_+^*\Or_{N^-}=p_-^*\Or_{N^-}
\end{equation}
If $E$ is a flat vector bundle over $M$ we get canonic
isomorphisms $p_-^*E|_\Sigma=i^*E=p_+^*E|_\Sigma$ using parallel 
transport along $X$.

The flow $\phi_t$ restricts to a smooth action of $\R$ on $\M$.
This action is free and we have a smooth orbit space
$\T:=\M/\R$. The mapping $p_+:\M\to\Sigma$
factors to a smooth map $\pi_+:\T\to\Sigma$. The fiber bundle
$p_-:\M\to\Sigma$ factors to a smooth fiber bundle
$\pi_-:\T\to\Sigma$. Let $q:\M\to\T$ denote the projection
and let $[X]\subseteq V_\M$ denote the line bundle over $\M$ generated
by the vector field $X$. We have a short exact sequence of vector bundles:
$$
0\to[X]\to V_\M\to q^*V_\T\to0
$$
So the rank of the vertical bundle $V_\T$ is
$\rank(V_\T)=\pi_-^*\ind-\pi_+^*\ind-1:\T\to\N_0$.
Since the line bundle $[X]$ is oriented via $X$ we obtain an isomorphism
of orientation bundles $\Or_{V_\M}=q^*\Or_{V_\T}$, that is we use the
$X$ last convention here. Using the isomorphism \eqref{E:AB}
%$\Or_{V_\M}\otimes p_+^*\Or_{N^-}=p_-^*\Or_{N^-}$
we get an isomorphism 
$q^*(\Or_{V_\T}\otimes\pi_+^*\Or_{N^-})=q^*\pi_-^*\Or_{N^-}$ which descends
to a canonic isomorphism of orientation bundles:
\begin{equation}\label{E:tt}
%o_\T:
\Or_{V_\T}\otimes\pi_+^*\Or_{N^-}=\pi_-^*\Or_{N^-}
\end{equation}
If $E$ is a flat vector bundle over $M$ then the isomorphism
$p_-^*E|_\Sigma=p_+^*E|_\Sigma$ descends to a canonic
isomorphism of flat vector bundles:
\begin{equation}\label{E:tttt}
\pi_-^*E|_\Sigma=\pi_+^*E|_\Sigma
\end{equation}

For two connected components $S$ and $S'$ of $\Sigma$ set
$\T(S,S'):=\pi_-^{-1}(S)\cap\pi_+^{-1}(S')\subseteq\T$.
Clearly this is the disjoint union over several connected 
components of $\T$. The rank of its  
vertical bundle $V_{\T(S,S')}$ is constant 
and $\rank(V_{\T(S,S')})=\ind(S)-\ind(S')-1$.
Note that this implies $\T(S,S')=\emptyset$ unless
$\ind(S)>\ind(S')$, a consequence of the strong transversality 
condition.
Also note, that we have $\T(S,S)=\emptyset$ by our very definitions.

\subsection{Compactification of ${\T}$}
\label{SS:comp_T}

Let $X$ be a Morse--Bott--Smale vector field on $M$.
For $k\geq0$ define the manifold of \emph{unparameterized 
$k$--times broken trajectories} by:
$$
\hat\T_k:=\underbrace{\T\times_\Sigma\T\times_\Sigma\cdots
\times_\Sigma\T}_{\text{$k+1$ copies}}
$$
Here 
$\T\times_\Sigma\T:=\{(x,y)\in\T\times\T\mid\pi_+(x)=\pi_-(y)\}$,
and similarly for more factors. Define a smooth map
$(\hat\pi_+)_k:\hat\T_k\to\Sigma$ by composing the projection
to the last factor with $\pi_+:\T\to\Sigma$.
Define a smooth fiber bundle $(\hat\pi_-)_k:\hat\T_k\to\Sigma$
by composing the projection to the first factor with
$\pi_-:\T\to\Sigma$. Note that by definition $\hat\T_0=\T$,
$(\hat\pi_+)_0=\pi_+$ and $(\hat\pi_-)_0=\pi_-$.

Let $V_{\hat\T_k}$ denote the vertical
bundle of $(\hat\pi_-)_k:\hat\T_k\to\Sigma$. 
Define an isomorphism of orientation bundles
\begin{equation}\label{E:pre_or_iso_k}
%o_{\hat\T_1}=
(-)^{p_1^*\rank(V_\T)}:\Or_{V_{\hat T_1}}\otimes(\hat\pi_+)_1^*\Or_{N^-}
=(\hat\pi_-)_1^*\Or_{N^-}
\end{equation}
as follows. Consider the commutative diagram 
$$
\xymatrix{
\hat\T_1=\T\times_\Sigma\T 
\ar[r]^-{p_2}
\ar[d]_{p_1}
& 
\T 
\ar[r]^-{\pi_+}
\ar[d]^{\pi_-}
&
\Sigma
\\
\T 
\ar[r]^{\pi_+}
\ar[d]_{\pi_-}
& 
\Sigma
\\
\Sigma
}
$$
where $p_1$ projects onto the 
first component and $p_2$ projects onto the second component of 
$\hat\T_1=\T\times_\Sigma\T$. 
Notice that $(\hat\pi_+)_1=\pi_+\circ p_2$ and
$(\hat\pi_-)_1=\pi_-\circ p_1$.
Clearly we have $p_2^*V_\T=V_{p_1}$ and thus 
$p_2^*\Or_{V_\T}=\Or_{V_{p_1}}$. Moreover the short exact sequence
$0\to V_{p_1}\to V_{\hat\T_1}\to p_1^*V_\T\to0$ provides
$\Or_{V_{p_1}}\otimes p_1^*\Or_{V_\T}
=\Or_{V_{\hat\T_1}}$. Finally recall the isomorphism \eqref{E:tt}.
%$o_\T:\Or_{V_{\pi_-}}\otimes\pi_+^*\Or_{N^-}=\pi_-^*\Or_{N^-}$.
Now define the isomorphism \eqref{E:pre_or_iso_k} as $(-)^{p_1^*\rank(V_\T)}$
times the following composition of isomorphisms:
\begin{eqnarray*}
\Or_{V_{\hat T_1}}\otimes(\hat\pi_+)_1^*\Or_{N^-}
&=&
\Or_{V_{p_1}}\otimes p_1^*\Or_{V_\T}\otimes p_2^*\pi_+^*\Or_{N^-}
\\&=&
p_2^*\Or_{V_\T}\otimes p_1^*\Or_{V_\T}\otimes 
p_2^*\pi_+^*\Or_{N^-}
\\&=&
p_1^*\Or_{V_\T}\otimes 
p_2^*(\Or_{V_\T}\otimes\pi_+^*\Or_{N^-})
\\&=&
p_1^*\Or_{V_\T}\otimes p_2^*\pi_-^*\Or_{N^-}
\\&=&
p_1^*(\Or_{V_\T}\otimes\pi_+^*\Or_{N^-})
\\&=&
p_1^*\pi_-^*\Or_{N^-}
\\&=&
(\hat\pi_-)_1^*\Or_{N^-}
\end{eqnarray*}

Suppose $E$ is a flat vector bundle over $M$.
In a similar manner, using \eqref{E:tttt}, we 
obtain a canonic isomorphism of flat vector bundles
\begin{equation}\label{E:ttt}
(\hat\pi_-)_1^*E|_\Sigma=(\hat\pi_+)_1^*E|_\Sigma
\end{equation}
as the following composition of isomorphisms:
$$
(\hat\pi_-)^*_1E|_\Sigma
=p_1^*\pi_-^*E|_\Sigma
=p_1^*\pi_+^*E|_\Sigma
=p_2^*\pi_-^*E|_\Sigma
=p_2^*\pi_+^*E|_\Sigma
=(\hat\pi_+)^*_1E|_\Sigma
$$

Set $\hat\T:=\bigsqcup_{k\geq0}\hat\T_k$. 
Define $\hat\pi_\pm:\hat\T\to\Sigma$ by $(\hat\pi_\pm)_k$ on $\hat\T_k$.
For connected components
$S,S'\subseteq\Sigma$ introduce a topology on 
$\hat\T(S,S'):=\hat\pi_-^{-1}(S)\cap\hat\pi_+^{-1}(S')$
by parameterizing the possibly broken trajectories with the help of $f$ and
taking the topology induced from $C^0([a,b],M)$, where
$a=f(S')$, $b=f(S)$. Topologize 
$\hat\T=\bigsqcup_{S,S'\subseteq\Sigma}\hat\T(S,S')$ 
as the disjoint union.

\begin{proposition}\label{P:comp_T}
The space $\hat\T$ is compact.
\end{proposition}

\begin{proof}
Clearly it suffices to show that $\hat\T(S,S')$ is compact,
for $\Sigma$ has only finitely many connected components.
The compactness of $\hat\T(S,S')$ follows from the theorem of
Arzela--Ascoli.
Indeed, the subset $\hat\T(S,S')\subseteq C^0([a,b],M)$ is
equicontinuous. For $t\in[a,b]$ which is non-critical this follows
from the fact that these are flow lines of a locally non-vanishing
vector field.
For critical $t$ one has to use the very explicit form of the 
singularities of $X$, \cf \cite{BH01}.

To see that $\hat\T(S,S')\subseteq C^0([a,b],M)$ is
closed suppose $\sigma_k\in\hat\T(S,S')$ converges uniformly to
$\sigma\in C^0([a,b],M)$. Then $\sigma$ too is parametrized by $f$, can only
pass through finitely many critical points, departs at $S$, ends at $S'$ and
away from the critical points $\sigma$ too is a flow line of $\frac X{|X|}$.
Hence $\sigma\in\hat\T(S,S')$.
\end{proof}

Recall that a smooth manifold with corners $P$ of dimension $n$
is a paracompact Hausdorff space modeled on 
$\{(x_1,\dotsc,x_n)|x_i\geq0\}\subseteq\R^n$.
The subset of points which in one (and hence every) chart
have exactly $k$ coordinates zero is called the \emph{$k$--corner
of $P$} and denoted by $\partial_kP$. It inherits the structure
of a smooth manifold of dimension $n-k$. The subset 
$\partial P:=\partial_1P\cup\cdots\cup\partial_nP$ is closed and
inherits the structure of a topological manifold.
The pair $(P,\partial P)$ is a topological manifold with boundary.

\begin{theorem}\label{T:comp_T}
Let $X$ be a Morse--Bott--Smale vector field on $M$.
There is a canonic way to equip the space
$\hat\T$ with the structure of a (compact) smooth manifold 
with corners such that the following hold:
\begin{enumerate}
\item\label{T:comp_Ti}
The inclusion $\hat\T_k\subseteq\hat\T$ is a diffeomorphism
onto the $k$--corner.
\item\label{T:comp_Tii}
The mapping $\hat\pi_+:\hat\T\to\Sigma$ is smooth.
\item\label{T:comp_Tiii}
The mapping $\hat\pi_-:\hat\T\to\Sigma$ is a smooth fiber 
bundle whose fibers are (compact) smooth manifolds with corners.
\item\label{T:comp_Tiv}
$\rank(V_{\hat\T})=\hat\pi_-^*\ind-\hat\pi_+^*\ind-1:\hat\T\to\N_0$, where
$V_{\hat\T}$ denotes the vertical bundle of $\hat\pi_-:\hat\T\to\Sigma$.
\item\label{T:comp_Tv}
If $X=-\grad_g(f)$ and $\hat\T(S,S')\neq\emptyset$ then $f(S)>f(S')$.
\item\label{T:comp_Tvi}
The isomorphism \eqref{E:tt} extends (uniquely) to a smooth isomorphism 
of orientation bundles
$\Or_{V_{\hat\T}}\otimes\hat\pi_+^*\Or_{N^-}
=\hat\pi_-^*\Or_{N^-}$.
Via $\Or_{V_{\hat\T}}|_{\hat\T_1}=\Or_{V_{\hat\T_1}}$
its restriction to $\hat\T_1$ identifies with \eqref{E:pre_or_iso_k}.
\item\label{T:comp_Tvii}
If $E$ is a flat vector bundle over $M$ then the isomorphism
\eqref{E:tttt} extends (uniquely) to a smooth isomorphism of flat vector
bundles $\hat\pi_-^*E|_\Sigma=\hat\pi_+^*E|_\Sigma$. Its restriction to
$\hat\T_1$ identifies with \eqref{E:ttt}.
\end{enumerate}
\end{theorem}

The proof of Theorem~\ref{T:comp_T} is contained in 
Appendix~\ref{A:proof_comp}.

\subsection{Compactification of $W^-$}\label{SS:comp_W}

For $k\geq0$ define a smooth manifold
$$
(\hat W^-)_k
:=\underbrace{\T\times_\Sigma\T\times_\Sigma\cdots
\times_\Sigma\T}_{\text{$k$ factors}}\times_\Sigma W^-
$$
Here $\T\times_\Sigma W^-
=\{(x,y)\in\T\times W^-\mid\pi_+(x)=p_-(y)\}$ and similarly
for more factors. 
Define a smooth mapping $\hat i_k:(\hat W^-)_k\to M$
by composing the projection to the last factor with 
$i_-:W^-\to M$. Define a smooth fiber bundle
$(\hat p_-)_k:(\hat W^-)_k\to\Sigma$ by composing the 
projection to the first factor with the mapping 
$\pi_-:\T\to\Sigma$ in the case $k\geq1$ and composing with
$p_-:W^-\to\Sigma$ if $k=0$. Note that by definition we have
$(\hat W^-)_0=W^-$, $\hat i_0=i_-$ and $(\hat p_-)_0=p_-$.

Let $V_{(\hat W^-)_k}$ denote the vertical bundle of
$(\hat p_-)_k:(\hat W^-)_k\to\Sigma$. Define an isomorphism
of orientation bundles
\begin{equation}\label{E:I}
%o_{(\hat W^-)_1}=
(-)^{p_1^*\rank(V_\T)}:
\Or_{V_{(\hat W^-)_1}}=(\hat p_-)_1^*\Or_{N^-}
\end{equation}
as follows. Consider the commutative diagram
\begin{equation}\label{D1}
\xymatrix{
(\hat W^-)_1=\T\times_\Sigma W^- 
\ar[r]^-{p_2} 
\ar[d]_{p_1}
&
W^- 
\ar[r]^-{i_-} 
\ar[d]^{p_-}
&
M
\\
\T 
\ar[r]^-{\pi_+}
\ar[d]_{\pi_-}
& 
\Sigma
\\
\Sigma
}
\end{equation}
where $p_1$ and $p_2$ denote the projection to the first and second
factor, respectively. Notice that $\hat i_1=i_-\circ p_2$ and
$(\hat p_-)_1=\pi_-\circ p_1$. Clearly we have $p_2^*V_{W^-}=V_{p_1}$ and
thus $p_2^*\Or_{V_{W^-}}=\Or_{V_{p_1}}$. Moreover the short exact sequence
$$
0\to V_{p_1}\to V_{(\hat W^-)_1}\to p_1^*V_\T\to0
$$
provides $\Or_{V_{p_1}}\otimes p_1^*\Or_{V_\T}=\Or_{V_{(\hat W^-)_1}}$.
Finally recall \eqref{E:A} and \eqref{E:tt}. Now define the isomorphism
\eqref{E:I} as $(-)^{p_1^*\rank(V_\T)}$ times the following composition
of isomorphisms:
\begin{eqnarray*}
\Or_{V_{(\hat W^-)_1}}
&=&
\Or_{V_{p_1}}\otimes p_1^*\Or_{V_\T}
\\&=&
p_2^*\Or_{V_{W^-}}\otimes p_1^*\Or_{V_\T}
\\&=&
p_2^*p_-^*\Or_{N^-}\otimes p_1^*\Or_{V_\T}
\\&=&
p_1^*(\pi_+^*\Or_{N^-}\otimes\Or_{V_\T})
\\&=&
p_1^*\pi_-^*\Or_{N^-}
\\&=&
(\hat p_-)^*_1\Or_{N^-}
\end{eqnarray*}

Suppose $E$ is a flat vector bundle over $M$. In a similar manner, using
\eqref{E:B} and \eqref{E:tttt}, define an isomorphism of flat vector 
bundles
\begin{equation}\label{E:II}
\hat i_1^*E=(\hat p_-)^*_1E|_\Sigma
\end{equation}
as the following composition of isomorphisms:
$$
\hat i^*_1E
=p_2^*i_-^*E
=p_2^*p_-^*E|_\Sigma
=p_1^*\pi_+^*E|_\Sigma
=p_1^*\pi_-^*E|_\Sigma
=(\hat p_-)^*_1E|_\Sigma
$$

Set $\hat W^-:=\bigsqcup_{k\geq0}(\hat W^-)_k$. 
Define $\hat i:\hat W^-\to M$ using $\hat i_k$ on $(\hat W^-)_k$.
Define $\hat p_-:\hat W^-\to\Sigma$ using $(\hat p_-)_k$ on
$(\hat W^-)_k$. 
Let $S\subseteq\Sigma$ be a connected component.
We are going to topologize $\hat W^-_S:=\hat p_-^{-1}(S)$.
Set $b:=f(S)$ and $a:=\min\{f(x)\mid x\in M\}$. Every element
$x\in\hat W^-_S$ can be considered as an unparametrized path $\sigma$ 
from $S$ to $\hat i(x)$ which we parametrize with the help of $f$,
\ie $\sigma\in C^0([a,b],M)$
$$
f(\sigma(t))=
\begin{cases}
a+b-t & \text{for $a\leq t\leq a+b-f(\hat i(x))$ and} 
\\
f(\hat i(x)) & \text{for $a+b-f(\hat i(x))\leq t\leq b$.}
\end{cases}
$$
We equip $\hat W^-_S$ with the topology induced from $C^0([a,b],M)$
and topologize $\hat W^-=\bigsqcup_{S\subseteq\Sigma}\hat W^-_S$
as the disjoint union.

\begin{proposition}\label{P:comp_W}
The space $\hat W^-$ is compact.
\end{proposition}

The proof of Proposition~\ref{P:comp_W} is exactly the same
as the proof of Proposition~\ref{P:comp_T}.

\begin{theorem}\label{T:comp_W}
Let $X$ be a Morse--Bott--Smale vector field on $M$.
There is a canonic way to equip the space $\hat W^-$ with the structure
of a (compact) smooth manifold with corners such that the following hold:
\begin{enumerate}
\item\label{T:comp_Wi}
The inclusion $(\hat W^-)_k\subseteq\hat W^-$ is a diffeomorphism
onto the $k$--corner of $\hat W^-$.
\item\label{T:comp_Wii}
The mapping $\hat i:\hat W^-\to M$ is smooth.
\item\label{T:comp_Wiii}
The mapping $\hat p_-:\hat W^-\to\Sigma$ is a smooth fiber
bundle whose fibers are (compact) smooth manifolds with corners.
\item\label{T:comp_Wiv}
$\rank(V_{\hat W^-})=\hat p_-^*\ind:\hat W^-\to\N_0$,
where $V_{\hat W^-}$ denotes the vertical bundle of 
$\hat p_-:\hat W^-\to\Sigma$.
\item\label{T:comp_Wv}
If $X=-\grad_g(f)$ then $\hat f(\hat W_S^-)\leq f(S)$ for all connected components
$S\subseteq\Sigma$, where $\hat f:=f\circ\hat i:\hat W^-\to\R$.
\item\label{T:comp_Wvi}
The isomorphism \eqref{E:A} extends (uniquely) to a smooth isomorphism
of orientation bundles
$\Or_{V_{\hat W^-}}=\hat p_-^*\Or_{N^-}$.
Via $\Or_{V_{\hat W^-}}|_{(\hat W^-)_1}=\Or_{V_{(\hat W^-)_1}}$ its
restriction to $(\hat W^-)_1$ identifies with \eqref{E:I}.
\item\label{T:comp_Wvii}
If $E$ is a flat vector bundle over $M$ then the isomorphism \eqref{E:B}
extends (uniquely) to a smooth isomorphism of flat vector bundles
$\hat i^*E=\hat p_-^*E|_\Sigma$. Its restriction to $(\hat W^-)_1$
identifies with \eqref{E:II}.
\end{enumerate}
\end{theorem}

The proof of Theorem~\ref{T:comp_W} is contained in 
Appendix~\ref{A:proof_comp}.

  \section{The geometric complex}
  \label{S:geom_complex}

\subsection{Fiber integration}\label{SS:fiber}

We will heavily use fiber integration. So let us fix 
a sign convention, which will be the same as in \cite{BT82}.
Suppose $\pi:F\to B$ is a smooth
fiber bundle. Here $F$, $B$ and the fibers $F_x=\pi^{-1}(x)$
may have corners and for simplicity we assume the fibers to 
be compact. However none of these manifolds is assumed to
be connected or of constant dimension.\footnote{Thus
even if we fiber integrate a form of homogeneous degree
the outcome need not be of homogeneous degree.}
Let $V_F\to F$ denote the vertical bundle of 
$\pi$. Finally suppose we have a flat vector bundle $E\to B$.
Given a differential form $\alpha\in\Omega(F;\pi^*E\otimes\Or_{V_F})$
we are going to define a form 
$\pi_*\alpha\in\Omega(B;E)$ as follows. Suppose $x\in B$
and $X_1,\dotsc,X_k\in T_xB$. Then the degree $k$ component
of $\pi_*\alpha\in\Omega(B;E)$ is defined by:
$$
(\pi_*\alpha)(X_1,\dotsc,X_k)
:=\int_{F_x}i_{\tilde X_k}\cdots i_{\tilde X_1}\alpha
$$
Here $\tilde X_i$ are vector fields on the total space of
$F$ defined along the fiber $F_x$ which project to $X_i$.
The pull back of $i_{\tilde X_k}\cdots i_{\tilde X_1}\alpha$
to the fiber $F_x$ is a form 
$\Omega(F_x;\pi^*E_x\otimes\Or_{F_x})$ which can be integrated 
over $F_x$ to yield an element of $E_x$. The outcome does not depend
on the choice of $\tilde X_i$.
For the fiber bundle $\pi:F\to B$ the fiber integration
$$
\pi_*:\Omega(F;\pi^*E\otimes\Or_{V_F})\to\Omega(B;E)
$$
has the following properties which
follow easily from our definition.

Suppose $f:B'\to B$ is a smooth mapping. Then we have
a pull back bundle $f^*\pi:f^*F\to B'$ and a smooth mapping
$\pi^*f:f^*F\to F$ which restricts to a diffeomorphism
$(f^*F)_x\to F_{f(x)}$ for every $x\in B'$. 
$$
\xymatrix{
f^*F
\ar[r]^-{\pi^*f}
\ar[d]_{f^*\pi}
&
F
\ar[d]^{\pi}
\\
B'
\ar[r]^-{f}
&
B
}
$$
Hence its tangent 
mapping induces an isomorphism of vector bundles
$V_{f^*F}=(\pi^*f)^*V_F$ over the total space of $f^*F$. 
Particularly we get a
canonic isomorphism of orientation bundles
$\Or_{f^*V_F}=(\pi^*f)^*\Or_{V_F}$ and thus a canonic
isomorphism of flat vector bundles
$(\pi^*f)^*(\pi^*E\otimes\Or_{V_F})
=(f^*\pi)^*f^*E\otimes\Or_{V_{f^*F}}$. Using this 
identification we have
$$
(f^*\pi)_*\circ(\pi^*f)^*=f^*\circ\pi_*:
\Omega(F;\pi^*E\otimes\Or_{V_F})\to
\Omega(B';f^*E)
$$
the \emph{naturality of fiber integration}.

Suppose $\tilde\pi:\tilde F\to F$ and $\pi:F\to B$ are two 
fiber bundles. Then
$\pi\circ\tilde\pi:\tilde F\to B$ is a fiber bundle too
and we have a canonic short exact sequence of vector bundles
$0\to V_{\tilde\pi}\to V_{\pi\circ\tilde\pi}\to\tilde\pi^*V_\pi\to0$
over $\tilde F$. This provides us with a canonic isomorphism
of orientation bundles 
$\Or_{V_{\pi\circ\tilde\pi}}
=\Or_{V_{\tilde\pi}}\otimes\tilde\pi^*\Or_{V_\pi}$ and thus with
a canonic isomorphism
$(\pi\circ\tilde\pi)^*E\otimes\Or_{V_{\pi\circ\tilde\pi}}
=\tilde\pi^*(\pi^*E\otimes\Or_{V_\pi})\otimes\Or_{V_{\tilde\pi}}$.
Via this identification we have
$$
\pi_*\circ\tilde\pi_*=(\pi\circ\tilde\pi)_*:
\Omega(\tilde F;(\pi\circ\tilde\pi)^*E\otimes\Or_{V_{\pi\circ\tilde\pi}})
\to
\Omega(B;E)
$$
the (covariant) \emph{functoriality of fiber integration}.

Suppose $\pi:F\to B$ is a fiber bundle and suppose
$E'$ and $E$ are two flat vector bundles over $B$.
For $\alpha\in\Omega(F;\pi^*E\otimes\Or_{V_F})$
and $\beta\in\Omega(B;E')$ we have
$$
\pi_*(\pi^*\beta\wedge\alpha)=\beta\wedge\pi_*\alpha
\in\Omega(B;E'\otimes E)
$$
a \emph{push--pull formula}. Here we only use the natural
isomorphism $\pi^*E'\otimes\pi^* E=\pi^*(E'\otimes E)$.

Suppose $F\to B$ is a fiber bundle and suppose $B$ has no
boundary. Then the $1$--corner 
$\pi|_{\partial_1F}:\partial_1F\to B$ is a smooth fiber bundle
with possibly non-compact fibers. Using the \emph{outward
pointing normal first} convention we get a canonic isomorphism of 
orientation
bundles $\Or_{V_F}|_{\partial_1F}=\Or_{V_{\partial_1F}}$ and
thus a canonic isomorphism of flat vector bundles
$(\pi^*E\otimes\Or_{V_F})|_{\partial_1F}
=(\pi|_{\partial_1F})^*E\otimes\Or_{V_{\partial_1F}}$.
Note that this is the same isomorphisms as the one we get from
the short exact sequence
$$
0\to V_{\partial_1F}\to V_F|_{\partial_1F}
\to V_F|_{\partial_1F}/V_{\partial_1F}\to0
$$
with the right hand side quotient oriented by the outward pointing normal.
Via this identification we find 
$$
\pi_*\circ d=d\circ\pi_*+A\circ(\pi|_{\partial_1F})_*:
\Omega(F;\pi^*E\otimes\Or_{V_F})\to\Omega(B;E)
$$
where $A:=(-)^k:\Omega^k(B;E)\to\Omega^k(B;E)$ is a grading
homomorphism.

\subsection{Integration and the differential}\label{SS:diff}

Let $X$ be a Morse--Bott--Smale vector field on a closed manifold $M$, 
\cf Definition~\ref{D:MBS}.
Let $E$ be a flat vector bundle over $M$. Define a vector space 
$$
C(X;E):=\Omega(\Sigma;E|_\Sigma\otimes\Or_{N^-})
$$
and an integration 
$$
\Int:=(\hat p_-)_*\hat i^*=(p_-)_*i_-^*:\Omega(M;E)\to C(X;E).
$$
Note that for $\alpha\in\Omega(M;E)$ we have
$\hat i^*\alpha\in\Omega(\hat W^-;\hat i^*E)$ which via the isomorphisms in
Theorem~\ref{T:comp_W}\itemref{T:comp_Wvi} and \itemref{T:comp_Wvii}
can be interpreted as an element of
$\Omega(\hat W^-;\hat p_-^*E|_\Sigma\otimes\hat p_-^*\Or_{N^-}\otimes\Or_{V_{\hat W^-}})$
and thus $(\hat p_-)_*\hat i^*\alpha$ is indeed an element of
$\Omega(\Sigma;E|_\Sigma\otimes\Or_{N^-})=C(X;E)$.
Similarly, using the isomorphisms in 
Theorem~\ref{T:comp_T}\itemref{T:comp_Tvi} and \itemref{T:comp_Tvii}, 
we define
$$
u:=(\hat\pi_-)_*\hat\pi_+^*=(\pi_-)_*\pi_+^*:C(X;E)\to C(X;E)
$$
and
$$
\delta:=d+uA:C(X;E)\to C(X;E)
$$
where $A:=(-)^k:\Omega^k(\Sigma;E|_\Sigma\otimes\Or_{N^-})\to
\Omega^k(\Sigma;E|_\Sigma\otimes\Or_{N^-})$ is a grading homomorphism, 
\cf section~\ref{SS:fiber}.

\begin{proposition}\label{P:int_homo}
$\Int\circ d=\delta\circ\Int$.
\end{proposition}

\begin{proof}
From section~\ref{SS:fiber} we obtain
$(\hat p_-)_*d=d(\hat p_-)_*+A((\hat p_-)|_{\partial_1\hat W^-})_*$ 
and thus:
$$
\Int\circ d=d\circ\Int+A((\hat p_-)|_{\partial_1\hat W^-})_*\hat i^*
$$
Using the description of the $1$--boundary of $\hat W^-$ from
Theorem~\ref{T:comp_W} and the obvious relation
$(\pi_-)_*(-)^{\rank(V_\T)}\pi_+^*=AuA$ we get:
\begin{eqnarray*}
((\hat p_-)|_{\partial_1\hat W^-})_*\hat i^*
&=&
(\pi_-\circ p_1)_*(-)^{p_1^*\rank(V_\T)}(i_-\circ p_2)^*
\\&=&
(\pi_-)_*(-)^{\rank(V_\T)}(p_1)_*p_2^*i_-^*
\\&=&
(\pi_-)_*(-)^{\rank(V_\T)}\pi_+^*(p_-)_*i_-^*
\\&=&
AuA\Int
\end{eqnarray*}
Therefore $\Int d=d\Int+A^2uA\Int=(d+uA)\Int=\delta\Int$.
\end{proof}

Since the integration $\Int:\Omega(M;E)\to C(X;E)$ is onto we must have
$\delta^2=0$ as a consequence of $d^2=0$. 
However, we will give a direct argument below.

\begin{proposition}
$\delta^2=0$.
\end{proposition}

\begin{proof}
Exactly as in the proof of Proposition~\ref{P:int_homo} we derive
$ud=du+uAu$ using the description of the $1$--boundary of $\hat\T$
from Theorem~\ref{T:comp_T}.
Then 
$\delta^2=duA+uAd+uAuA=(du-ud+uAu)A=0$, where we also made use of the
obvious $dA=-Ad$.
\end{proof}

\begin{definition}[Geometric complex]
The complex $(C(X;E),\delta)$ is called the \emph{geometric complex
with values in the flat vector bundle $E$}.
It is a generalization of the deRham complex and the classical
Morse complex at the same time. The homomorphism of chain complexes
$\Int:\Omega(M;E)\to C(X;E)$ is called the \emph{integration homomorphism}.
\end{definition}

\subsection{Grading and filtration}\label{SS:filt}

The geometric complex $C(X;E)$
can be graded and filtered in a compatible way. 
For a connected component $S\subseteq\Sigma$ let us write
$E_S:=E_\Sigma|_S=(E|_\Sigma\otimes\Or_{N^-})|_S$. Define
$$
C_p^q(X;E):=\bigoplus_{\ind(S)\geq p}\Omega^{q-\ind(S)}(S;E_S)
$$
where the direct sum is over all connected components $S\subseteq\Sigma$ of
index larger or equal to $p$. Let us also introduce the notation
$C^q(X;E):=C^q_0(X;E)$, and note that $\bigoplus_qC^q(X;E)=C(X;E)$.
Finally let us define $C_p(X;E)=\bigoplus_qC^q_p(X;E)$.
Clearly $C^q_p(X;E)=C^q(X;E)\cap C_p(X;E)$.
Of course the grading and the filtration stem from a bigrading on $C(X;E)$.
However, since the bigrading will not be compatible with the differential
we will not consider it in the sequel.

\begin{proposition}\label{P:filt_d}
The differential
of the geometric complex preserves the grading and the 
filtrations as indicated below:
$$
\delta:C^q_p(X;E)\to C^{q+1}_p(X;E)
$$
\end{proposition}

\begin{proof}
This is obviously true for the deRham differential. So it suffices to
consider $u=(\pi_-)_*(\pi_+)^*$. Let $\alpha'$ be a form
of degree $|\alpha'|$ supported on a connected component $S'\subseteq\Sigma$. Let
$S\subseteq\Sigma$ be another connected component and let $\alpha$
denote the restriction of $(\pi_-)_*(\pi_+)^*\alpha'$ to $S$.
That is we integrate $(\pi_+)^*\alpha'$ only over $\T(S,S')$.
We may assume $\T(S,S')\neq\emptyset$.
Since $\rank(V_{\T(S,S')})=\ind(S)-\ind(S')-1$, 
see Theorem~\ref{T:comp_T}\itemref{T:comp_Tiv}, 
the degree of $\alpha$ will be $|\alpha'|-\ind(S)+\ind(S')+1$.
This shows $\delta:C^q\to C^{q+1}$. Since it is the rank of a vector bundle
we must have $\ind(S)-\ind(S')-1\geq0$ and
thus $\delta:C_p\to C_p$. 
\end{proof}

A Morse--Bott--Smale vector field on $M$ also provides a filtration of
$\Omega^*(M)$ as follows. For a connected component $S\subseteq\Sigma$ 
let $\hat W^-_S$ denote its compactified unstable manifold and let
$\hat i_S:\hat W^-_S\to M$ denote the restriction of the smooth map
$\hat i:\hat W^-\to M$, \cf Theorem~\ref{T:comp_W}. Let
$\Omega^q_p(M;E)$ denote the space of forms
$\alpha\in\Omega^q(M;E)$ which vanish on some open neighborhood of
$M(p):=\bigcup_{\ind(S)<p}\hat i_S(\hat W^-_S)$. Clearly
$\Omega^*_{p+1}(M;E)\subseteq\Omega^*_p(M;E)$ and this
filtration is preserved by the deRham differential.
Note that $(\Omega^*_p(M;E),d)$ computes the relative cohomology
$H^*(M,M(p);E)$.

\begin{proposition}\label{P:filt_I}
The integration preserves the grading and the filtration
as indicated below.
$$
\Int:\Omega^q_p(M;E)\to C^q_p(X;E)
$$
\end{proposition}

\begin{proof}
Let $\alpha'$ be a form on $M$ of homogeneous degree $|\alpha'|$ and
let $S\subseteq\Sigma$ be a connected component. Let $\alpha$ denote
the restriction of $(p_-)_*i^*\alpha'$ to $S$. That is we only integrate
over $W^-_S$. Since
$\rank(V_{W_S^-})=\ind(S)$, see Theorem~\ref{T:comp_W}\itemref{T:comp_Wiv} 
we must have $|\alpha|=|\alpha'|-\ind(S)$.
This implies $\Int:\Omega^q\to C^q$. The rest is obvious.
\end{proof}

\subsection{Cohomology of the geometric complex}\label{SS:thom}

Let us consider the filtered graded complex 
$C^q_p(X;E)$, see Proposition~\ref{P:filt_d}. It gives rise to a spectral sequence
$(E_kC^q_p(X;E),\delta_k)$. Since $C^*_0(X;E)=C^*(X;E)$ and $C^*_{\dim M+1}(X;E)=0$
the spectral sequence will converge to $HC^*(X;E)$. Moreover, for the
$E_1$--term we obviously have $E_1C^q_p(X;E)=H^{q-p}(\Sigma_p;E_{\Sigma_p})$.
Here $\Sigma_p=\bigsqcup_{\ind(S)=p}S$ denotes the critical 
manifold of index $p$ and
$E_{\Sigma_p}=E_\Sigma|_{\Sigma_p}=(E|_\Sigma\otimes\Or_{N^-})|_{\Sigma_p}$.
This proves the second part of Theorem~\ref{T:thom}.

Via this identification the differential on the $E_1$--term has a particularly 
simple form. All components but
\begin{equation}\label{E:donE1}
\delta_1:H^{q-p}(\Sigma_p;E_{\Sigma_p})\to H^{q-p}(\Sigma_{p+1};E_{\Sigma_{p+1}})
\end{equation}
vanish, and \eqref{E:donE1} is given by
$\delta_1=(-)^{q-p}(\hat\pi_-)_*\circ(\hat\pi_+)^*$. Here 
$$
(\hat\pi_+)^*:H^*(\Sigma_p;E_{\Sigma_p})\to
H^*\bigl(\hat\T(\Sigma_{p+1},\Sigma_p);(\hat\pi_+)^*E\bigr)
$$ 
is given by pull back along $\hat\pi_+:\hat\T(\Sigma_{p+1},\Sigma_p)\to\Sigma_p$, and
$$
(\hat\pi_-)_*:H^*\bigl(\hat\T(\Sigma_{p+1},\Sigma_p);(\hat\pi_+)^*E\bigr)\to
H^*(\Sigma_{p+1};E_{\Sigma_{p+1}})
$$ 
is the map induced by fiber integration along
$\hat\pi_-:\hat\T(\Sigma_{p+1},\Sigma_p)\to\Sigma_{p+1}$.
For this note that since the index difference
is one, the fibers of the bundle
$\hat\pi_-:\hat\T(\Sigma_{p+1},\Sigma_p)\to\Sigma_{p+1}$ are zero dimensional,
hence there is no boundary, and fiber integration does indeed induce a
mapping in cohomology.

\begin{remark}
The differentials $\delta_k$ for $k\geq1$ should be considered as 
generalization of the classical \emph{counting of instantons}.
These are differentials of finite dimensional complexes and can
be considered the \emph{combinatorics of the dynamic}.
If $X$ is a Morse--Smale vector field then $(E_1C^*,\delta_1)$ is the 
familiar Morse complex and $\delta_k=0$ for all $k\geq2$.
\end{remark}

Recall the filtration $\Omega^q_p(M;E)$ introduced in section~\ref{SS:filt}.
Clearly we have $\Omega^*_0(M;E)=\Omega^*(M;E)$ and 
$\Omega^*_{\dim M+1}(M;E)=0$. So the associated spectral sequence
$(E_k\Omega^*(M;E),\delta_k)$ will converge to $H^*(M;E)$. 
Let $W^-_{\Sigma_p}$ denote the unstable manifold of the critical manifold 
of index $p$. Using the description of the boundary of $\hat W^-_{\Sigma_p}$ in
Theorem~\ref{T:comp_W} and the fact that $i:W^-_{\Sigma_p}\to M$ is an
embedding we obtain a homomorphism
\begin{equation}\label{E:ext}
\Omega^*_p(M;E)/\Omega^*_{p+1}(M;E)\xrightarrow{i^*}
\Omega^*_c(W^-_{\Sigma_p};i^*E)
\end{equation}
where the latter denotes forms with compact support.
This is a quasi isomorphism, \ie induces an isomorphism in cohomology, 
since $\Omega^*_p(M;E)$ computes $H^*(M,M(p);E)$.
It follows from the Thom isomorphism theorem that $\Int:E_1\Omega^q_p(M;E)\to
E_1C^q_p(X;E)$, \cf Proposition~\ref{P:filt_I},
is an isomorphism. The first part of Theorem~\ref{T:thom} now follows from a
standard spectral sequence argument.

\begin{remark}\label{R:specOC}
Note that this argument actually shows that 
$$
\Int:(E_k\Omega^*_p(M;E),\delta_k)\to(E_kC^*_p(X;E),\delta_k)
$$ 
is an isomorphism for all $k\geq1$. That is,
the spectral sequences induced by $\Omega^*_p(M;E)$ and
$C^*_p(X;E)$ are isomorphic.
\end{remark}

\begin{corollary}\label{C:euler}
For the Euler characteristics we have:
\begin{equation}\label{E:euler}
\chi(M)
=\sum_q(-)^q\chi(\Sigma_q)
%=\sum_{S\subseteq\Sigma}(-)^{\ind(S)}\chi(S)
\end{equation}
\end{corollary}

\begin{proof}
The Euler characteristics $\chi(M)$ equals the Euler characteristics of
$H^*(M)$ which equals the Euler characteristics of the geometric complex
$\chi(HC^*)$ by Theorem~\ref{T:thom}\itemref{T:thom:i}. Since the spectral sequence
$E_kC^*$ converges to $HC^*$ we get $\chi(HC^*)=\chi(E_1C^*)$ which implies
\eqref{E:euler} in view of Theorem~\ref{T:thom}\itemref{T:thom:ii}.
\end{proof}

\begin{remark}
Suppose $X$ is a vector field whose zero set $\Sigma$ is a closed submanifold
but not necessarily non-degenerate. For every connected component
$S\subseteq\Sigma$ one can define an index $\IND(S)=\IND_X(S)$ generalizing the
Hopf index, see \cite{D95}, and one has $\chi(M)=\sum_{S\subseteq\Sigma}\IND(S)$.
If the vector field $-X$ is Morse--Bott then the indices are related by
$\IND_X(S)=(-)^{\ind_{-X}(S)}\chi(S)$. 
This also implies Corollary~\ref{C:euler} without using the Smale condition.
Finally notice that if the zero set of $X$ is degenerate
$\IND(S)$ may be non-zero even if $\chi(S)=0$.

The generalized Hopf index can easily be defined as follows. Let $U$ be a
connected open neighborhood of $S$. Then the vector field provides a mapping
$X:(U,U\setminus S)\to(TU,TU\setminus U)$. The Hopf index then is the unique
integer which makes the following diagram commutative:
$$
\xymatrix{
H^n(TU,TU\setminus U;\pi^*_{TU}\Or_U) 
\ar[r]^-{X^*}
\ar@{=}[d]_{\textrm{Thom}}
&
H^n(U,U\setminus S;\Or_U)
\ar@{=}[d]^{\textrm{Thom}}
\\
H^0(U)=\R
\ar[r]^-{\IND(S)}
&
\R=H^{\dim(S)}(S;\Or_S)
}
$$
\end{remark}

\begin{corollary}[Morse inequalities]
For every $q_0\in\Z$ we have
$$
\sum_{q\geq q_0}(-)^{q-q_0}b^q(M)
\leq
\sum_{q\geq q_0}(-)^{q-q_0}\sum_pb^{q-p}(\Sigma_p)
$$
where $b^q(M):=\dim H^q(M;E)$ and
$b^q(\Sigma_p)=\dim H^q(\Sigma_p;E_{\Sigma_p})$.
\end{corollary}

\begin{proof}
Recall that for a finite dimensional complex $C^*$ and any 
$q_0\in\Z$ we have:
$$
\sum_{q\geq q_0}(-)^{q-q_0}\dim HC^q\leq
\sum_{q\geq q_0}(-)^{q-q_0}\dim C^q
$$
Applying this to the complex $(E_kC^*(X;E),\delta_k)$
we obtain
$$
\sum_{q\geq q_0}(-)^{q-q_0}\dim E_{k+1}C^q\leq
\sum_{q\geq q_0}(-)^{q-q_0}\dim E_kC^q
$$
which in turn implies:
$$
\sum_{q\geq q_0}(-)^{q-q_0}\dim E_\infty C^q\leq
\sum_{q\geq q_0}(-)^{q-q_0}\dim E_1C^q
$$
The statement now follows from Theorem~\ref{T:thom}.
\end{proof}

  \section{Torsion}
  \label{S:torsion}

\subsection{Euler structures}\label{SS:euler}

The importance of Euler 
structures for torsion comes from the fact that they permit to remove the 
metric ambiguity of analytic and geometric torsion, see below.
In the case $\chi(M)=0$ Euler structures where introduced by Turaev 
\cite{Tu90} for exactly this purpose. A detailed discussion of the 
general case can be found in \cite{BH}. 
In this paper we will only consider Euler structures over $\R$. Hence
the space of Euler structures will be an affine version of $H_1(M;\R)$,
\ie $H_1(M;\R)$ acts free and transitively on the set of Euler structures.

In \cite{BH} the set of Euler structures $\Eul_{x_0}(M)$ based at $x_0\in M$
was defined as the set of equivalence classes $[X,c]$, where
$X$ is a vector field with non-degenerate zeros and $c\in C_1(M;\R)$ is a
singular one chain with $\partial c=\sum_{x\in\mathcal X}\IND(x)x-\chi(M)x_0$. 
Here $\mathcal X$ denotes the (finite) set of zeros of $X$ and 
$\IND(x)=\IND_X(x)$ its
Hopf index at $x$. Notice that $\IND(x)=\pm1$ for we assume $X$ to be
non-degenerate. Every such chain $c$ is called \emph{Euler chain} or
\emph{Turaev spider} for $X$ based at $x_0$. Since $\sum_{x\in\mathcal
X}\IND(x)=\chi(M)$ every such vector field admits Euler chains.
By definition $[X_1,c_1]=[X_2,c_2]$ iff 
$c_2-c_1=c(X_1,X_2)$ modulo boundaries. Here
$c(X_1,X_2)$ is the singular one chain obtained from the zero set of
a non-degenerate homotopy from $X_1$ to $X_2$. Note that
$$
\partial c(X_1,X_2)
=\sum_{x\in\mathcal X_2}\IND_{X_2}(x)x
-\sum_{x\in\mathcal X_1}\IND_{X_1}(x)x
$$
and that changing the homotopy 
changes $c(X_1,X_2)$ by a boundary only. The action of
$[\sigma]\in H_1(M;\R)$ on $[X,c]\in\Eul_{x_0}(M)$ is defined by
$[X,c]+[\sigma]:=[X,c+\sigma]$. This action clearly is free and transitive.
If $x_1$ is another base point and $\tau$ a path from $x_0$ to $x_1$
then 
\begin{equation}\label{E:tor_base}
[X,c]\mapsto[X,c]-\chi(M)\tau:=[X,c-\chi(M)\tau]
\end{equation}
defines 
an $H_1(M;\R)$--equivariant 
(affine) bijection $\Eul_{x_0}(M)\to\Eul_{x_1}(M)$. Clearly this does
only depend on the homology class of the path from $x_0$ to $x_1$.
If $\chi(M)=0$ it does not depend on the base point at all and the set
of Euler structures is defined without reference to a base point, \cf
\cite{Tu90}.

We will also use another way to represent Euler structures which is 
better suited to remove ambiguities from
the analytic torsion. For a base point $x_0\in M$ define $\Eul^*_{x_0}(M)$
as the set of equivalence classes $[g,\alpha]$, where $g$ is a Riemannian
metric on $M$ and $\alpha\in\Omega^{\dim M-1}(M\setminus\{x_0\};\Or_M)$ with
$d\alpha=E_g$, where $E_g\in\Omega^{\dim M}(M;\Or_M)$ denotes the Euler form
of $g$. By definition $[g_1,\alpha_1]=[g_2,\alpha_2]$ iff
$\alpha_2-\alpha_1=\cs(g_1,g_2)$ modulo $d\bigl(\Omega^{\dim
M-1}(M;\Or_M)\bigr)$, where $\cs(g_1,g_2)$ denotes the Chern--Simon invariant.
There is a free and transitive action of $H^{\dim M-1}(M;\Or_M)$ on 
$\Eul^*_{x_0}$ defined by $[g,\alpha]+[\beta]:=[g,\alpha-\beta]$. Recall that
Poincar\'e duality provides an isomorphism 
$H_1(M;\R)=H^{\dim M-1}(M;\Or_M)$, where $[\sigma]\in H_1(M;\R)$ corresponds 
to $[\beta]\in H^{\dim M-1}(M;\Or_M)$
iff $\int_\sigma\omega=\int_M\omega\wedge\beta$ for all closed one forms
$\omega\in\Omega^1(M;\R)$. There is a canonic affine isomorphism
$\Eul_{x_0}(M)=\Eul^*_{x_0}(M)$ which can be regarded as affine version of 
Poincar\'e duality. To define it we need two preliminary propositions.
For details see \cite{BH}.

For a closed one form $\omega\in\Omega^1(M;\R)$, a Riemannian metric $g$,
a form $\alpha\in\Omega^{\dim M-1}(M\setminus\{x_0\};\Or_M)$ with
$d\alpha=E_g$ on $M\setminus\{x_0\}$ and a smooth function $f:M\to\R$ such 
that $\omega':=\omega-df$ vanishes in a neighborhood of $x_0$ define:
$$
\cS(\omega,g,\alpha;f)
:=\int_{M\setminus\{x_0\}}\omega'\wedge\alpha-\int_MfE_g+\chi(M)f(x_0)
$$

\begin{proposition}\label{P:S}
The quantity $\cS(\omega,g,\alpha;f)$ does not depend on the choice of $f$
and will be denoted by $\cS(\omega,g,\alpha)$. It should be considered 
a regularization of the possibly divergent integral
$\int_{M\setminus\{x_0\}}\omega\wedge\alpha$.
It has the following properties:
\begin{enumerate}
\item\label{S:i}
$\cS(\omega,g,\alpha)$ is linear in $\omega$. 
\item\label{S:ii}
$\cS(dh,g,\alpha)=-\int_MhE_g+\chi(M)h(x_0)$ for all smooth functions
$h:M\to\R$.
\item\label{S:iv}
$\cS(\omega,g,\alpha+\beta)-\cS(\omega,g,\alpha)=\int_M\omega\wedge\beta$
for all closed $\beta\in\Omega^{\dim M-1}(M;\Or_M)$.
\item\label{S:iii}
If $[\alpha_1,g_1]=[\alpha_2,g_2]\in\Eul_{x_0}^*(M)$ then
$\cS(\omega,g_2,\alpha_2)-\cS(\omega,g_1,\alpha_1)
=\int_M\omega\wedge\cs(g_1,g_2)$ for all closed one forms $\omega$.
\end{enumerate}
\end{proposition}

\begin{proof}
The independence of $f$ follows from a straight forward calculation
using $d\alpha=E_g$, $\int_ME_g=\chi(M)$ and Stokes' theorem.
All of the properties are obvious, except \itemref{S:iii}. For this one use
$d\cs(g_1,g_2)=E_{g_2}-E_{g_1}$ and Stokes' theorem.
For details see \cite[Lemma~2]{BH}. 
\end{proof}

For a closed one form $\omega\in\Omega^1(M;\R)$, a vector field $X$ with
non-degenerate zero set $\mathcal X$, a Riemannian metric $g$ and a smooth 
function $f:M\to\R$ such that $\omega':=\omega-df$ vanishes in a 
neighborhood of $\mathcal X$ define:
$$
\cR(\omega,X,g;f):=
\int_{M\setminus\mathcal X}\omega'\wedge X^*\Psi_g
-\int_MfE_g+\sum_{x\in\mathcal X}\IND(x)f(x)
$$
Here $\Psi_g$ denotes the global angular form, \cf \cite{BT82}, also called
Mathai--Quillen form in \cite{BZ92}. Recall that $\Psi_g$ is define as follows:
Let $n=\dim(M)$ and let $\pi_{TM}:TM\to M$ be the tangent bundle equipped with 
the Levi--Civita connection. Let  $\Vol_g\in\Omega^n(TM;\pi^*_{TM}\Or_M)$, the global
volume form which vanishes when contracted with horizontal vectors on $TM$ and which
assigns to an $n$-tuple of vertical vectors `their volume times their orientation'.
Moreover let $\xi$ denote the Euler vector field on $TM$ which
assigns to a point $v\in TM$ the vertical vector $-v\in T_v(TM)$. Then
$$
\Psi_g
:=\frac{\Gamma(n/2)}{(2\pi)^{n/2}|\xi|^n}
i_\xi\Vol_g\in\Omega^{n-1}(TM\setminus M;\pi^*_{TM}\Or_M).
$$

\begin{proposition}\label{P:R}
The quantity $\cR(\omega,X,g;f)$ does not depend on the choice of $f$
and will be denoted by $\cR(\omega,X,g)$.\footnote{The quantity
$\cR(\omega,X,g)$ was considered by Bismut--Zhang \cite{BZ92} for
$X=-\grad_g(h)$ and $h$ a Morse function.}
It should be considered a regularization of the possibly divergent integral
$\int_{M\setminus\mathcal X}\omega\wedge X^*\Psi_g$. It has the following
properties:
\begin{enumerate}
\item\label{R:i}
$\cR(\omega,X,g)$ is linear in $\omega$.
\item\label{R:ii}
$\cR(dh,X,g)=-\int_MhE_g
+\sum_{x\in\mathcal X}\IND(x)h(x)$
for all smooth functions $h:M\to\R$.
\item\label{R:iii}
$\cR(\omega,X,g_2)-\cR(\omega,X,g_1)=\int_M\omega\wedge\cs(g_1,g_2)$
for two Riemannian metrics $g_1$ and $g_2$.
\item\label{R:iv}
$\cR(\omega,X_2,g)-\cR(\omega,X_1,g)=\int_{c(X_1,X_2)}\omega$
for two vector fields $X_1$ and $X_2$ with non-degenerate zeros.
\end{enumerate}
\end{proposition}

\begin{proof}
The independence of $f$ follows from straight forward calculation using
$d\Psi_g=E_g$, Stokes' theorem and 
$\lim_{\epsilon\to0}\int_{\partial B_\epsilon(x)}X^*\Psi_g=-\IND(x)$
with $B_\epsilon(x)$ denoting the $\epsilon$--ball around $x$.
The first two properties are obvious. For \itemref{R:iii} use
$\Psi_{g_2}-\Psi_{g_1}=\cs(g_1,g_2)$, $E_{g_2}-E_{g_1}=d\cs(g_1,g_2)$
and Stokes' theorem. Concerning \eqref{R:iv} it is easy to see that 
it suffices to consider the case where there is a homotopy
$\mathbb X$ from $X_1$ to $X_2$ such that $\omega$ is exact on a
neighborhood of the zero set of $\mathbb X$. Let $p:I\times M\to M$ denote
the projection. The homotopy $\mathbb X$ provides a section of $p^*TM$
which we may assume to be transversal to the zero section. Choose $f$ such
that $p^*\omega'=p^*\omega-p^*df$ vanishes on the zero set of $\mathbb X$.
Then 
$$
\sum_{x\in\mathcal X_2}\IND_{X_2}(x)f(x)
-\sum_{x\in\mathcal X_1}\IND_{X_1}(x)f(x)
=\int_{c(X_1,X_2)}df=\int_{c(X_1,X_2)}\omega
$$
and
\begin{eqnarray*}
\int_M\omega'\wedge X_2^*\Psi_g-\int_M\omega'\wedge X_1^*\Psi_g
&=&
\int_{I\times M}d\bigl(p^*\omega'\wedge\mathbb X^*\tilde p^*\Psi_g\bigr)
\\&=&
-\int_{I\times M}p^*(\omega'\wedge E_g)=0.
\end{eqnarray*}
Here $\tilde p:p^*TM\to TM$ denotes the canonic vector bundle homomorphism
over $p:I\times M\to M$. More details can be found in \cite[Section~3]{BH}.
\end{proof}

One can now define the affine isomorphism $\Eul_{x_0}(M)=\Eul_{x_0}^*(M)$
by saying that $[X,c]$ corresponds to $[g,\alpha]$ iff
$$
\cR(\omega,X,g)-\cS(\omega,\alpha,g)-\int_c\omega=0
$$
for all closed one forms $\omega\in\Omega^1(M;\R)$. In view of the last two
Propositions this is indeed a well defined isomorphism affine over the 
Poincar\'e duality $H_1(M;\R)=H^{\dim M-1}(M;\Or_M)$.
We will from now on drop the notation $\Eul^*_{x_0}(M)$ and feel free to
represent an Euler structure $\e\in\Eul_{x_0}(M)$ as an equivalence class 
$[X,c]$ or as an equivalence class $[g,\alpha]$.

For the Morse--Bott case we introduce yet another way to
represent Euler structures. This can be considered a localization
of the Euler structure at the critical manifold of a Morse--Bott
vector field. Fix a Morse--Bott vector field $X$
with critical manifold $\Sigma$ and define $\Eul^X_{x_0}(M)$ as the set of 
equivalence classes $[X,c,\{\e_S\}]$, where
for every connected component $S\subseteq\Sigma$, 
$\e_S\in\Eul_{x_S}(S)$ with $x_S\in S$ and $c\in C_1(M;\R)$ is a one chain with
$\partial c=\sum_{S\subseteq\Sigma}\IND(S)x_S-\chi(M)x_0$. Since we have
$\sum_{S\subseteq\Sigma}\IND(S)=\chi(M)$ such $c$ indeed exist.
Given two such data $(X,c^1,\{\e_S^1\})$ and
$(X,c^2,\{\e_S^2\})$ with $\e_S^i\in\Eul_{x_S^i}(S)$
we choose paths $\tau_S$ in $S$ from
$x_S^1$ to $x_S^2$. Then $\e_S^1-\chi(S)\tau_S\in\Eul_{x_S^2}(S)$
and hence $(\e_S^1-\chi(S)\tau_S)-\e_S^2\in H_1(S;\R)$.
We define $[X,c^1,\{e_S^1\}]=[X,c^2,\{e_S^2\}]$
if for one (and hence every) choice of paths $\tau_S$ as above we have
$$
\sum_{S\subseteq\Sigma}(-)^{\ind_{-X}(S)}\bigl((\e_S^1-\chi(S)\tau_S)-\e_S^2\bigr)
=c^2-c^1-\sum_{S\subseteq\Sigma}\IND(S)\tau_S
$$
in $H_1(M;\R)$.
This is indeed independent of the paths $\{\tau_S\}$ for we have
$\IND(S)=(-)^{\ind_{-X}(S)}\chi(S)$ for every 
$S\subseteq\Sigma$. For $[\sigma]\in H_1(M;\R)$ the action is defined by
$[X,c,\{\e_S\}]+[\sigma]:=[X,c+\sigma,\{\e_S\}]$.
Note that for $[\sigma_S]\in H_1(S;\R)$ we have
\begin{equation}\label{E:sss}
[X,c,\{\e_S+[\sigma_S]\}]
=[X,c+\sum_{S\subseteq\Sigma}(-)^{\ind_{-X}(S)}\sigma_S,\{\e_S\}]
\end{equation}

Note how this specializes to the previous definition 
if $X$ has non-degenerate zeros only, for in this case there is a unique 
base point and a unique Euler structure on every $S=\{\pt\}$.
Also note how this localization degenerates for the Morse--Bott--Smale
vector field $X\equiv0$.

Suppose $[X,c,\{\e_S\}]$ as above and suppose the Euler structures
$\e_S\in\Eul_{x_S}(S)$ are represented as $\e_S=[Y_S,c_S]$. With the help of a tubular
neighborhood of $\Sigma$, a connection and a fiber metric on the normal bundle
of $\Sigma$ and a bump function we can construct a globally defined vector field 
$Y$ which restricts to $Y_S$ and is supported in a small neighborhood of 
$\Sigma$. Consider the homotopy $t\mapsto X+tY$. Using the Morse--Bott
assumption for $X$ we see that for sufficiently small
$\epsilon>0$ the zeros $\mathcal X^+$ of $X^+:=X+\epsilon Y$ will all be non-degenerate and
coincide with the set of zeros of $Y_S$. For the index of such a 
zero $x$ we have $\ind^M_{-X^+}(x)=\ind^S_{-Y_S}(x)+\ind_{-X}^M(S)$. Therefor 
$c^+:=c+\sum_{S\subseteq\Sigma}(-)^{\ind_{-X}(S)}c_S$ will be an Euler chain for $X^+$,
and hence $[X^+,c^+]\in\Eul_{x_0}(M)$. This Euler structure does not depend
on the choices we made and the construction provides a natural affine 
isomorphism $\Eul^X_{x_0}(M)=\Eul_{x_0}(M)$. Obviously we have

\begin{lemma}\label{L:euler}
Suppose $X$ is a Morse--Bott vector field,
$\e\in\Eul_{x_0}(M)$ and
$\e_S\in\Eul_{x_S}(S)$ for every connected component $S\subseteq\Sigma$.
Then there exists an Euler chain $c\in C_1(M;\R)$ such that $\e=[X,c,\{\e_S\}]$.
Moreover, every other such Euler chain differs from $c$ by a boundary.
\end{lemma}

For a Morse--Bott vector field $X$, a closed one form
$\omega\in\Omega^1(M;\R)$, a Riemannian metric $g$ and a Morse vector field
$Y$ on $\Sigma$ define
$$
\cR(\omega,X,g;Y):=\cR(\omega,X^+,g)
-\sum_{S\subseteq\Sigma}(-)^{\ind_{-X}(S)}\cR(\omega,Y|_S,g|_S)
$$
where $X^+=X+\epsilon Y$ is a small perturbation of $X$ as above.

\begin{proposition}\label{P:RR}
The quantity $\cR(\omega,X,g;Y)$ does not depend on $Y$ nor on the other choices
entering the definition of the perturbation $X^+$ and will be denoted
by $\cR(\omega,X,g)$. It should be
considered a regularization of the possibly divergent integral 
$\int_{M\setminus\Sigma}\omega\wedge X^*\Psi_g$. 
It has the following properties:
\begin{enumerate}
\item\label{RR:i}
$\cR(\omega,X,g)$ is linear in $\omega$.
\item\label{RR:ii}
$\cR(dh,X,g)=-\int_MhE_g
+\sum_{S\subseteq\Sigma}(-)^{\ind_{-X}(S)}\int_ShE_{g|_S}$
for all smooth functions $h:M\to\R$.
\item\label{RR:iii}
$\cR(\omega,X,g_2)-\cR(\omega,X,g_1)=\int_M\omega\wedge\cs(g_1,g_2)
-\sum_{S\subseteq\Sigma}(-)^{\ind_{-X}(S)}\int_S\omega\wedge\cs(g_1|_S,g_2|_S)$
for two Riemannian metrics $g_1$ and $g_2$.
\item\label{RR:iv}
If there exists a function $f$ such that $\omega'=\omega-df$ vanishes in a
neighborhood of $\Sigma$ then:
$$
\cR(\omega,X,g)=\int_{M\setminus\Sigma}\omega'\wedge X^*\Psi_g-\int_MfE_g
+\sum_{S\subseteq\Sigma}(-)^{\ind_{-X}(S)}\int_SfE_{g|_S}
$$
\end{enumerate}
\end{proposition}

\begin{proof}
All this is an immediate consequences of Proposition~\ref{P:R}.
For the independence of $Y$ note that a homotopy 
from $Y_1$ to $Y_2$ as vector fields on $\Sigma$ provides a homotopy
from $X_1^+$ to $X_2^+$ via a similar construction. The zero set of this
homotopy will coincide with the zero set of the homotopy from $Y_1$ to $Y_2$
we started with. However, the orientations
they inherit differ by $(-)^{\ind_{-X}(S)}$.
Now apply Proposition~\ref{P:R}\itemref{R:iv} twice.
\end{proof}

\begin{proposition}\label{P:equal_euler}
Suppose $\e'=[g,\alpha]$ and
$\e=[X,c,\{\e_S\}]$
are two Euler structures based at $x_0$ and suppose
$\e_S=[g|_S,\alpha_S]$.
Then $\e=\e'$ iff
$$
0=
\cR(\omega,X,g)
-\cS(\omega,\alpha,g)
-\int_c\omega
+\sum_{S\subseteq\Sigma}(-)^{\ind_{-X}(S)}
\cS(\omega,\alpha_S,g|_S)
$$
for all closed one forms $\omega\in\Omega^1(M;\R)$.
\end{proposition}

\begin{proof}
This is immediate from our definitions.
\end{proof}

\subsection{Analytic torsion}\label{SS:an_tor}

Let $M$ be a closed manifold, $E$  a flat vector bundle over $M$ and let
$\e\in\Eul_{x_0}(M)$ be a an Euler structure based at $x_0\in M$. 
Consider the deRham complex $\Omega^*(M;E)$, the determinant line
of its cohomology $\det H^*(M;E)$ and 
the one dimensional vector space:
\begin{equation}\label{E:DetME}
\Det_{x_0}(M;E):=\det H^*(M;E)\otimes(\det E_{x_0})^{-\chi(M)}
\end{equation}
The analytic torsion will be a non-zero element in $\Det_{x_0}(M;E)$
defined up to multiplication by $\pm1$, equivalently a scalar product
on $\Det_{x_0}(M;E)$, associated with the Euler structure $\e$.

To define the analytic torsion choose a Riemannian metric $g$ on $M$, 
a fiber metric $\mu$ on $E$ and
$\alpha\in\Omega^{\dim M-1}(M\setminus\{x_0\};\Or_M)$ so that $d\alpha=E_g$
on $M\setminus\{x_0\}$ and so that $\e=[g,\alpha]\in\Eul_{x_0}(M)$, \cf
section~\ref{SS:euler}.

The deRham complex $\Omega(M;E)$ inherits a scalar product which permits to
define the formal adjoint $d^\sharp:\Omega(M;E)\to\Omega(M;E)$
of the deRham differential $d$ and the Laplacian
$\Delta:=d\circ d^\sharp+d^\sharp\circ d$. Clearly
$\Delta^q:\Omega^q(M;E)\to\Omega^q(M;E)$ is a symmetric
positive semi-definite elliptic second order differential operator.
Therefor we have the \emph{Ray--Singer torsion}, see \cite{RS71}:
$$
\log T_\an^{E,g,\mu}:=\frac12\sum_q(-)^{q+1}q\log\det{}'\Delta^q\in\R 
$$
The fiber metric $\mu$ on $E$ determines a fiber metric on its
determinant bundle $\det E$ and a canonic length one section $e$ of
$\det E\otimes\Or_E$. There is a unique one form
$\theta=\theta_\mu^E\in\Omega^1(M;\R)$ such that 
$\nabla^{\det E\otimes\Or_E}_Ye=\theta(Y)e$.\footnote{The orientation
bundle $\Or_E$ just comes in to ensure the existence of a canonic globally
defined length one section $e$. The bundle $\det E$ has locally defined
length one sections with a sign ambiguity. One could as well use these
to first locally define $\theta$. Then one has to observe that the sign
ambiguity does not affect $\theta$ and that the locally defined one forms
actually define a global one form.
To circumvent these inconveniences we made use of $\Or_E$.}
Since the connection on $\det E\otimes\Or_E$ is flat $\theta$ is a 
closed one form. It measures to what extend the fiber metric $\mu$ 
on $\det E$ is not parallel. Indeed, one can easily check
$\theta(Y)=-\frac12\tr_\mu\nabla_Y^E\mu$, which could serve as alternative
definition, \cf \cite{BH}, \cite{BZ92}. The cohomology class 
$\Theta_E:=[\theta]\in H^1(M;\R)$ does not depend on $\mu$ and 
is a complete obstruction to the existence of a fiber metric on $E$
such that the induced fiber metric on $\det E$ is parallel.\footnote{Note
that the theta in \cite{BZ92} is $-2$ times the $\theta$ here.}
For later use let us state

\begin{lemma}\label{L:XYZ}
\
\begin{enumerate}
\item\label{L:XYZi}
Suppose $c$ is a closed path starting and ending at $x_0$. Then
parallel transport along $c$ provides an automorphism of $\det E_{x_0}$
which is given by multiplication with $\pm e^{-\langle\Theta_E,[c]\rangle}$
where $[c]\in H_1(M;\R)$ denotes the homology class represented by $c$.
\item\label{L:XYZii}
Let $c$ be a path from $x_0$ to $x_1$, suppose $\mu$ is a fiber metric on
$E$ and let $||\cdot||_{\det E_{x_i}}^\mu$ denote the induced scalar 
product on $\det E_{x_i}$. Parallel transport along $c$ will map 
$||\cdot||_{\det E_{x_0}}^\mu$ to 
$e^{\int_c\theta_E^\mu}\cdot||\cdot||_{\det E_{x_1}}^\mu$.
\end{enumerate}
\end{lemma}

Define a scalar product on $H^*(M;E)$ by restricting 
the scalar product on $\Omega^*(M;E)$ to harmonic forms. Let
$||\cdot||_\Hodge^{E,g,\mu}$ denote the induced scalar product on
$\det H^*(M;E)$. Let $||\cdot||^\mu_{x_0}$
denote the scalar product on $(\det E_{x_0})^{\chi(M)}$ induced from
$\mu|_{E_{x_0}}$. Define a linear isomorphism 
$$
\tau_\an^{E,g,\alpha,\mu}:(\det E_{x_0})^{\chi(M)}=\det H^*(M;E)
$$
so that it maps $||\cdot||^\mu_{x_0}$ to
$T_\an^{E,g,\mu}e^{-\cS(\theta,g,\alpha)}||\cdot||_\Hodge^{E,g,\mu}$.
This is well defined up to sign and will be considered as element in
$\Det_{x_0}(M;E)$, \cf \eqref{E:DetME}.

\begin{proposition}\label{P:an_indep}
$\tau_\an^{E,g,\alpha,\mu}\in\Det_{x_0}(M;E)$
is independent of $\mu$ and does not depend on the 
representative of the Euler structure $\e=[g,\alpha]$. 
Hence we will denote it by $\tau_\an^{E,\e}$.
\end{proposition}

\begin{proof}
Let $g_1$ and $g_2$ be two Riemannian metrics on $M$, let
$\mu_1$ and $\mu_2$ be two fiber metrics on $E$ and suppose
$[g_1,\alpha_1]=\e=[g_2,\alpha_2]$, see section~\ref{SS:euler}.
Let $h$ denote the unique smooth function so that
$||\cdot||^{\mu_2}_{\det E}=e^h||\cdot||^{\mu_1}_{\det E}$ 
where $||\cdot||_{\det E}^{\mu_i}$ denotes the fiber metric on $\det E$
induced by $\mu_i$. Then we clearly have
$\log\bigl(||\cdot||_{x_0}^{\mu_2}/||\cdot||^{\mu_1}_{x_0}\bigr)=\chi(M)h(x_0)$.
Moreover, $-dh=\theta_2-\theta_1$ where $\theta_i:=\theta_E^{\mu_i}$.
With our notation the anomaly formula for the Ray--Singer torsion, see
\cite[Theorem~0.1]{BZ92} reads:
$$
\log\frac{T_\an^{E,g_2,\mu_2}||\cdot||_\Hodge^{E,g_2,\mu_2}}
{T_\an^{E,g_1,\mu_1}||\cdot||_\Hodge^{E,g_1,\mu_1}}
=\int_MhE_{g_1}+\int_M\theta_2\wedge\cs(g_1,g_2)
$$
Using Proposition~\ref{P:S} this implies
$$
\log\frac{e^{-\cS(\theta_2,g_2,\alpha_2)}T_\an^{E,g_2,\mu_2}||\cdot||_\Hodge^{E,g_2,\mu_2}}
{e^{-\cS(\theta_1,g_1,\alpha_1)}T_\an^{E,g_1,\mu_1}||\cdot||_\Hodge^{E,g_1,\mu_1}}
=\chi(M)h(x_0)
$$
and the proposition follows.
\end{proof}

\begin{definition}[Analytic torsion]\label{D:an_tor}
The element $\tau_\an^{E,\e}\in\Det_{x_0}(M;E)$ is well defined up to a 
sign and will be called the \emph{analytic torsion}.\footnote{
The sign ambiguity of the analytic torsion can be removed using a homology
orientation of $M$, \ie an orientation of $\det H^*(M;\R)$, see \cite{Tu02}.}
\end{definition}

\begin{remark}\label{R:an_tor_dep}
For every $\sigma\in H_1(M;\R)$ and every $\e\in\Eul_{x_0}(M)$ we have
$\tau_\an^{E,\e+\sigma}=\tau_\an^{E,\e}\cdot e^{-\langle\Theta_E,\sigma\rangle}$.
A path $\tau$ from $x_0$ to $x_1$ provides
an isomorphism $\Eul_{x_0}(M)=\Eul_{x_1}(M)$ 
denoted by $\e\mapsto\e-\chi(M)\tau$, see \eqref{E:tor_base} in section~\ref{SS:euler}. 
Moreover parallel transport along $\tau$ provides an identification 
$\Det_{x_0}(M;E)=\Det_{x_1}(M;E)$. Via the latter we have
$\tau^{E,\e}_\an=\tau_\an^{E,\e-\chi(M)\tau}$.
This follows from Proposition~\ref{P:S} and
Lemma~\ref{L:XYZ}\itemref{L:XYZii}.
\end{remark}

\subsection{Geometric torsion}\label{SS:geom_tor}

Let $X$ be Morse--Bott--Smale vector field on a closed manifold $M$,
let $E$ be a flat vector bundle over $M$ and let $\e\in\Eul_{x_0}(M)$ be an
Euler structure based at $x_0\in M$. Consider the geometric complex
$C^*(M;E)$, \cf section~\ref{SS:diff}, let $\det HC^*(X;E)$ denote the
determinant line of its cohomology and consider the one dimensional vector
space:
$$
\Det_{x_0}(X;E)=\det HC^*(X;E)\otimes(\det E_{x_0})^{-\chi(M)}
$$
The geometric torsion will be a non-zero element in
$\Det_{x_0}(X;E)$ defined up to multiplication by $\pm1$, equivalently a 
scalar product on $\Det_{x_0}(X;E)$, associated to the Euler
structure $\e$.

In order to define the geometric torsion choose an Euler structure
$\e_S\in\Eul_{x_S}(S)$ for every connected component $S\subseteq\Sigma$
and choose $c$ such that $\e=[-X,c,\{\e_S\}]$, see Lemma~\ref{L:euler}. 
This gives rise to an isomorphism, well defined up to multiplication by
$\pm1$,
\begin{equation}\label{E:base_pts}
(\det E_{x_0})^{\chi(M)}
\overset c=\bigotimes_{S\subseteq\Sigma}
\Bigl(\bigl(
\det(E_S)_{x_S}\bigr)^{\chi(S)}
\Bigr)^{(-)^{\ind(S)}}
\end{equation}
as follows. We choose a path $c_S$ from $x_0$ to $x_S$ and use parallel 
transport along $c_S$ to get an isomorphism $\det E_{x_0}=\det E_{x_S}$. 
Taking tensor products (and duals) and using Corollary~\ref{C:euler} we 
obtain an isomorphism
\begin{equation}\label{E:base_pts_pre}
(\det E_{x_0})^{\chi(M)}
\overset{\{c_S\}}\longrightarrow
\bigotimes_{S\subseteq\Sigma}
\bigl((\det E_{x_S})^{\chi(S)}\bigr)^{(-)^{\ind(S)}}
\end{equation}
Choose an orientation of $(\Or_{N^-})_{x_S}$ to get an
identification $(\Or_{N^-})_{x_S}=\R$ and hence an identification
$(E_S)_{x_S}=E_{x_S}$. 
Set 
$$
c_0
:=c-\sum_{S\subseteq\Sigma}(-)^{\ind(S)}\chi(S)c_S
\in C_1(M;\R)
$$
and note that $\partial c_0=0$. Now define \eqref{E:base_pts} as
$e^{-\langle\Theta_E,[c_0]\rangle}$ times the isomorphism 
\eqref{E:base_pts_pre}. It follows immediately from
Lemma~\ref{L:XYZ}\itemref{L:XYZi} that the resulting isomorphism
\eqref{E:base_pts} does not depend on the choice of paths $\{c_S\}$ and actually
only depends on $\e$ and $\{\e_S\}$ --- up to sign.

As explained in 
section~\ref{SS:filt} the geometric complex is a filtered graded complex 
$C^q_p(X;E)$. 
In view of Theorem~\ref{T:thom}\itemref{T:thom:ii} the associated spectral sequence 
converges to $HC^*(X;E)$. This provides us with a natural identification
$$
\det HC^*
=\det GHC^*
=\det E_\infty C^*
=\cdots
=\det E_2C^*
=\det E_1C^*
$$
where $C^*=C^*(X;E)$ and
$GHC^*$ denotes the associated graded of $HC^*$ which inherits a
filtration as usual from the filtration on $C^*$. For more details
on determinant lines see appendix~\ref{A:homalg}.
Together with the description of the $E_1$--term in
Theorem~\ref{T:thom}\itemref{T:thom:ii} this yields an isomorphism:
\begin{equation}\label{E:17}
\det HC^*(X;E)=\bigotimes_{S\subseteq\Sigma}\bigl(\det H^*(S;E_S)\bigr)^{(-)^{\ind(S)}}
\end{equation}
The tensor product of the analytic torsions provides an isomorphism
\begin{equation}\label{taus}
\bigotimes_{S\subseteq\Sigma}
\bigl((\det E_{x_S})^{\chi(S)}\bigr)^{(-)^{\ind(S)}}
=
\bigotimes_{S\subseteq\Sigma}\bigl(\det H^*(S;E_S)\bigr)^{(-)^{\ind(S)}}
\end{equation}
Let $\tau_\geom^{E,X,c,\{\e_S\}}\in\Det_{x_0}(X;E)$ denote the element
corresponding to the composition of \eqref{E:base_pts}, \eqref{taus} 
and \eqref{E:17}.
Using Remark~\ref{R:an_tor_dep} and \eqref{E:sss} we immediately obtain

\begin{proposition}\label{P:geom_tor_indep}
The element $\tau^{E,X,c,\{\e_S\}}\in\Det_{x_0}(X;E)$ does
not depend on the choice of $\e_S$, $c$ and hence will be denoted by
$\tau_\geom^{E,\e,X}$.
\end{proposition}

\begin{definition}[Geometric torsion]
The element $\tau_\geom^{E,\e,X}\in\Det_{x_0}(X;E)$
is well defined up to a sign
and called the \emph{geometric torsion}.\footnote{The sign ambiguity of the
geometric torsion can be removed using a homology orientation, that is an
orientation of $\det H^*(M;\R)$, see \cite{Tu02}.}
\end{definition}

\begin{remark}\label{R:geom_tor_dep}
For every $\sigma\in H_1(M;\R)$ and every $\e\in\Eul_{x_0}(M)$ we have
$\tau_\geom^{E,\e+\sigma,X}=\tau_\geom^{E,\e,X}\cdot e^{-\langle\Theta_E,\sigma\rangle}$.
A path $\tau$ from $x_0$ to $x_1$ provides
an isomorphism $\Eul_{x_0}(M)=\Eul_{x_1}(M)$ 
denoted by $\e\mapsto\e-\chi(M)\tau$, see \eqref{E:tor_base} in section~\ref{SS:euler}. 
Moreover parallel transport along $\tau$ provides an identification 
$\Det_{x_0}(X;E)=\Det_{x_1}(X;E)$. Via the latter we have
$\tau^{E,\e,X}_\geom=\tau_\geom^{E,\e-\chi(M)\tau,X}$.
Note that this is analogous to the properties of the analytic torsion, \cf
Remark~\ref{R:an_tor_dep}.
\end{remark}

\begin{remark}\label{R:MB_geom_tor}
Note how the geometric torsion specializes to the analytic torsion
for the Morse--Bott--Smale vector field $X\equiv0$.

Also note how the geometric torsion specializes to the combinatorial torsion 
for a Morse--Smale vector field $X$. 
More precisely, if $\e=[-X,c]$ then the geometric torsion coincides with 
the composition:
$$
(\det E_{x_0})^{\chi(M)}\overset c=\det C^*(X;E)=\det HC^*(X;E)
$$
This follows from Proposition~\ref{P:spec0} in appendix~\ref{A:homalg}, 
which is a very week statement in this special case.
\end{remark}

Recall from Theorem~\ref{T:thom} that the integration induces an 
isomorphism $H^*(M;E)=HC^*(X;E)$. So we get an induced isomorphism:
\begin{equation}\label{E:DetI}
\Det_{x_0}(M;E)=\Det_{x_0}(X;E)
\end{equation}
The following is the main result of this paper.

\begin{theorem}\label{T:main}
Let $X$ be a Morse--Bott--Smale vector field on a closed manifold $M$,
let $E$ a flat vector bundle over $M$ and suppose $\e$ is an Euler 
structure based at $x_0\in M$. 
Then the isomorphism \eqref{E:DetI} maps the analytic torsion
$\tau_\an^{E,\e}$ to the geometric torsion $\tau_\geom^{E,\e,X}$,
up to sign.
\end{theorem}

If $X$ is a Morse--Smale pair then this is due to
Bismut--Zhang, see \cite{BZ92}, which itself is an extension
of a result of Cheeger \cite{Ch77}, \cite{Ch79} and M\"uller \cite{Mu78}.
In this formulation it can be found in \cite[Theorem 5]{BH}.
The proof of Theorem~\ref{T:main} is contained in section~\ref{S:thm} below
and will make use of the classical Bismut--Zhang theorem, that is we will
assume Theorem~\ref{T:main} proved for Morse--Smale vector fields.

\subsection{Derivation of Theorem~\ref{T:intro} from Theorem~\ref{T:main}}

We continue to use the notation from section~\ref{SS:geom_tor}.
So we have an Euler structure $\e\in\Eul_{x_0}(M)$ and
Euler structures $\e_S\in\Eul_{x_S}(S)$ for every connected component
$S\subseteq\Sigma$ and we have a one chain $c$ such that $\e=[-X,c,\{\e_S\}]$.

Choose a fiber metric $\mu$ on $E$ and let $\mu_S$ denote the induced
fiber metric on $E_S$. Choose a Riemannian metric $g$ on $M$ and let
$g_S$ denote the restriction to $S\subseteq\Sigma$.
Choose $\alpha$ so that $\e=[g,\alpha]$ and for every connected component 
$S\subseteq\Sigma$ choose $\alpha_S$ with $\e_S=[g_S,\alpha_S]$, \cf
section~\ref{SS:euler}. From Proposition~\ref{P:equal_euler} we have:
\begin{equation}\label{E:AAA}
\cR(\theta,-X,g)
-\cS(\theta,\alpha,g)
-\int_c\theta
+\sum_{S\subseteq\Sigma}(-)^{\ind(S)}\cS(\theta,\alpha_S,g_S)
=0
\end{equation}
In view of Theorem~\ref{T:main} we have a commutative diagram:
%%%%%%%%%%%%%%%%%%%%%%%%%%%%%%%%%%%%%%%%%%%%%%%%
%$$
%\begin{CD}
%\bigotimes\bigl((\det E_S)_{x_S})^{\chi(S)}\bigr)^{(-)^{\ind(S)}}
%@<c<< 
%(\det E_{x_0})^{\chi(M)}
%@>\tau_\an^{E,\e}>>
%\det H^*(M;E)
%\\
%@V\otimes\tau_\an^{E_S,\e_S}VV @. @VV{\Int}V
%\\
%\bigotimes\bigl(\det H^*(S;E_S)\bigr)^{(-)^{\ind(S)}}
%@<<<
%\eqref{E:17}
%@<<<
%\det HC^*(X;E)
%\end{CD}
%$$
%%%%%%%%%%%%%%%%%%%%%%%%%%%%%%%%%%%%%%%%%%%%%%%
$$
\xymatrix{
\bigotimes\bigl((\det(E_S)_{x_S})^{\chi(S)}\bigr)^\pm
\ar[d]_-{\bigotimes\tau_\an^{E_S,\e_S}}
& 
(\det E_{x_0})^{\chi(M)}
\ar[l]_-{c}
\ar@/^1pc/[dr]^-{\tau_\an^{E,\e}}
\ar[d]^{\tau_\geom^{E,\e,X}}
\\
\bigotimes\bigl(\det H^*(S;E_S)\bigr)^\pm
&
\det HC^*(X;E)
\ar[l]_-{\eqref{E:17}}
&
\det H^*(M;E)
\ar[l]_-{\Int}
}
$$
The tensor products are over all connected components $S\subseteq\Sigma$
and the powers $\pm$ are abbreviations for $(-)^{\ind(S)}$.

Note that $\Int$ maps $||\cdot||_\Hodge^{E,g,\mu}$
to $(T_\met^M(E,X,g,\mu))^{-1}||\cdot||_\geom^{E,X,g,\mu}$ by the very
definition of the metric torsion. Applying
Proposition~\ref{P:lap_tor} in appendix~\ref{A:homalg}
once for every $E_k$--term we see that \eqref{E:17} maps
$||\cdot||_\geom^{E,X,g,\mu}$ to:
$$
(T_\comb^M(E,X,g,\mu))^{-1}\cdot
\bigotimes_{S\subseteq\Sigma}
\bigl(||\cdot||_\Hodge^{E_S,g_S,\mu_S}\bigr)^{(-)^{\ind(S)}}
$$
It follows from Lemma~\ref{L:XYZ}\itemref{L:XYZii} that $c$ maps 
$||\cdot||_{x_0}^\mu$ to:
$$
e^{\int_c\theta}\cdot\bigotimes\bigl(||\cdot||_{x_S}^{\mu_S}\bigr)^{(-)^{\ind(S)}}
$$
Therefore the commutativity of the diagram implies:
\begin{multline*}
T_\an^M(E,g,\mu)e^{-\cS(\theta,\alpha,g)}
(T_\met^M(E,X,g,\mu))^{-1}(T_\comb^M(E,X,g,\mu))^{-1}
\\
=e^{\int_c\theta}\cdot\prod_{S\subseteq\Sigma}\Bigl(
T_\an^M(E_S,g_S,\mu_S)
e^{-\cS(\theta,\alpha_S,g_S)}
\Bigr)^{(-)^{\ind(S)}}
\end{multline*}
Together with \eqref{E:AAA} Theorem~\ref{T:intro} follows at once.

\begin{remark}
The proof actually shows that Theorem~\ref{T:main} is equivalent to
Theorem~\ref{T:intro}. In view of the Bismut--Zhang theorem,
see \cite[Theorem~0.2]{BZ92} we may assume Theorem~\ref{T:main} proved
for Morse--Smale vector fields, see also \cite[Theorem 5]{BH}.
\end{remark}

  \section{Proof of main Theorem~\ref{T:main}}
  \label{S:thm}

\subsection{Homotopies}\label{SS:homotopy}

A Morse--Bott homotopy $\mathbb X=\{X_s\}_{s\in\R}$ from a 
Morse--Bott--Smale vector field $X^-$ to another Morse--Bott--Smale 
vector field $X^+$ is a smooth family of vector fields for which there 
exists $\epsilon>0$ so that $X_s=X^-$ when $s\leq-1+\epsilon$ and 
$X_s=X^+$ when $s\geq1-\epsilon$. To such a homotopy one associate a 
Morse--Bott vector field $Y$ on the smooth manifold 
$N:=M\times[-1,1]$ defined by
$$
Y(x,s):=X_s(x)+1/2(s^2-1)\frac{\partial}{\partial s}.
$$
The homotopy is called Morse--Bott--Smale homotopy if $Y$ is a 
Morse--Bott--Smale vector field. Be aware that it is in general
not possible to find a Morse--Bott--Smale homotopy $\X$ connecting 
two given Morse--Bott--Smale vector fields $X^-$ and $X^+$.

It is obvious from the definition of $Y$ that the discussion 
in sections~\ref{S:comp} and \ref{S:geom_complex} 
as well as the compacity of $\hat\T$ and $\hat W^-$ continue to hold.  
We will decorate $\Sigma$, $\T$, $u$ and $\Int$ by $X^\pm$ or 
$Y$ to indicate the relevant vector field.

We will use several canonic identifications. First, $\Sigma_Y$ 
identifies with the disjoint union of 
$\Sigma_{X^-}$ and $\Sigma_{X^+}$. Note that 
$N^-_Y|_{\Sigma_{X^-}}=N^-_{X^-}$ hence
$\ind_Y|_{\Sigma_{X^-}}=\ind_{X^-}:\Sigma_{X^-}\to\N_0$ and
$\Or_{N^-_Y}|_{\Sigma_{X^-}}=\Or_{N^-_{X^-}}$.
Moreover $N^-_{X^+}\subseteq N^-_Y|_{\Sigma_{X^+}}$ with a canonically 
oriented one dimensional quotient. Therefore 
$\ind_Y|_{\Sigma_{X^+}}=1+\ind_{X^+}:\Sigma_{X^+}\to\N_0$ and 
$\Or_{N^-_Y}|_{\Sigma_{X^-}}=\Or_{N^-_{X^+}}$. For the last one we use
the \emph{inward pointing normal last} convention.

Let $p:N\to M$ denote the
canonical projection and suppose $E$ is a flat vector bundle over $M$.
From our identifications we obtain
$$
C^*(Y;p^*E)=C^*(X^-;E)\oplus C^{*-1}(X^+;E)
$$ 
as graded vector spaces and
$\delta_Y=
\left(\begin{smallmatrix}
\delta_{X^-}&0
\\
v_\X&\delta_{X^+}
\end{smallmatrix}\right)$.
Define a grading homomorphism $B=(-)^k:C^k(X^+;E)\to C^k(X^+;E)$ and
let 
$$
\varphi_\X:=B\circ v_\X:C^*(X^-;E)\to C^*(X^+;E).
$$
The equation $\delta_Y^2=0$ implies 
$v_\X\circ\delta_{X^-}+\delta_{X^+}\circ v_\X=0$ and hence
\begin{equation}\label{E:phi_Y}
\varphi_\X\circ\delta_{X^-}=\delta_{X^+}\circ\varphi_\X.
\end{equation}

Moreover the integration $\Int_Y:\Omega(N;p^*E)\to C(Y;p^*E)$ decomposes as
$\Int_Y=
\left(\begin{smallmatrix}
\Int_{X^-}\circ\inc_-^*
\\
h_\X
\end{smallmatrix}\right)
$
where $\inc_\pm:M\to M\times\{\pm1\}\subseteq N$ denotes the canonical
inclusion. Slightly generalizing Proposition~\ref{P:int_homo} to the
boundary case we have:
\begin{equation}\label{E:int_Y}
\Int_Y\circ d=\delta_Y\circ\Int_Y-B\circ\Int_{X^+}\circ\inc_+^*
\end{equation}
The extra term stems from the additional $1$--boundary component of
$Y$'s unstable manifold which coincides with the main stratum
of the unstable manifold of $X^+$.
The sign $B$ stems partially from Stokes' theorem for
fiber integration, see section~\ref{SS:fiber}, and partially from
our identification $\Or_{N^-_Y}|_{\Sigma_{X^+}}=\Or_{N^-_{X^+}}$ 
using the inward pointing normal last convention. Define:
$$
H_\X:=B\circ h_\X\circ p^*:\Omega^*(M;E)\to C^{*-1}(X^+;E)
$$
Equation \eqref{E:int_Y} implies:
$$
h_\X\circ d
=v_\X\circ\Int_{X^-}\circ\inc_-^*
+\delta_{X^+}\circ h_\X
-B\circ\Int_{X^+}\circ\inc_+^*
$$
and thus:
\begin{equation}\label{E:H_Y}
H_\X\circ d+\delta_{X^+}\circ H_\X=\varphi_\X\circ\Int_{X^-}-\Int_{X^+}
\end{equation}

So we get an induced mapping in cohomology 
$\varphi_\X:HC^*(X^-;E)\to HC^*(X^+;E)$ which satisfies
$\varphi_\X\circ\Int_{X^-}=\Int_{X^+}$ and in
view of Theorem~\ref{T:thom}\itemref{T:thom:i} is an isomorphism
which does not depend on the particular homotopy from $X^-$ to $X^+$.
Notice also, that for a constant homotopy $\X$ we have
$\varphi_\X=\id$ on chain level.

Also note that $\varphi_\X:C^*_p(X^-;E)\to C^*_p(X^+;E)$ preserves the
filtration and hence it 
induces homomorphism intertwining the spectral sequences:
$\varphi_\X:E_kC^*(X^-;E)\to E_kC^*(X^+;E)$.

\subsection{A small perturbation}\label{SS:small_pert}

Let $X$ be a Morse--Bott--Smale vector field on a closed
manifold $M$ with critical manifold $\Sigma:=\Sigma_X$. 
Let $Y$ be a Morse--Smale vector field on $\Sigma$.
As in section~\ref{SS:euler} use a tubular neighborhood of $\Sigma$ 
to extend $Y$ to a globally defined vector field which is supported
in a small neighborhood of $\Sigma$. Define a homotopy $\X$ 
from $X$ to $X^+:=X+Y$ by $\X_t:=X+tY$.
By choosing the extension of $Y$ to $M$ carefully we may assume that 
the following hold true:
\begin{enumerate}
\item\label{BM:i}
$X^+$ is Morse--Smale and $\Sigma_{X^+}\subseteq\Sigma$.
\item\label{BM:ia}
We have $N^-_{X^+}=N^-_X\oplus N^-_{Y|_\Sigma}$,
$\ind_{X^+}^M=\ind^M_X|_{\Sigma_{X^+}}+\ind^\Sigma_{Y|_\Sigma}$ and
$\Or_{N^-_{X^+}}=\Or_{N^-_X}\otimes\Or_{N^-_{Y|_\Sigma}}$.
\item\label{BM:ii}
$X^+$ is tangential to $\Sigma$. For $x$ and $x'$ in a connected component
$S\subseteq\Sigma$ we have $\hat\T_{X^+}(x,x')=\hat\T_{Y|_S}(x,x')$.
\item\label{BM:iii}
Suppose $S\neq S'\subseteq\Sigma$ are two different connected components of
$\Sigma$, $x\in\Sigma_{X^+}\cap S$, $x'\in\Sigma_{X^+}\cap S'$ and suppose
$\hat\T_X(S,S')=\emptyset$. Then $\hat\T_{X^+}(x,x')=\emptyset$. 
\item\label{BM:iv}
The homotopy $\X$ is Morse--Bott--Smale.
\item\label{BM:v}
The homotopy $\X$ is tangential to $\Sigma$, \ie every $\X_t$ is tangential
to $\Sigma$. For every connected component $S\subseteq\Sigma$ and $x\in S$ we
have $\hat\T_\X(x,S)=\hat W_{Y|_S}^-(x)$.
\item\label{BM:vi}
Suppose $S\neq S'\subseteq\Sigma$ are two different connected components of
$\Sigma$, $x\in\Sigma_{X^+}\cap S$ and $\hat\T_X(S,S')=\emptyset$.
Then $\hat\T_\X(x,S')=\emptyset$.
\end{enumerate}

With the help of \itemref{BM:i} we can define an interesting filtration 
on the geometric complex $C(X^+;E)$ by:
$$
C(X^+;E)_p=\bigoplus_{p'\geq p}
\Omega(\Sigma_{X^+}\cap\Sigma_{p'};E\otimes\Or_{N^-_{X^+}})
$$
Because of \itemref{BM:iii} the differential on $C(X^+;E)$ preserves this
filtration. We thus get a spectral sequence $E_kC(X^+;E)$ converging to
$HC(X^+;E)$. In view of \itemref{BM:ia} and \itemref{BM:ii} we can 
compute the $E_1$--term as:
$$
E_1C^q(X^+;E)_p=
HC^{q-p}(Y|_{\Sigma_p};E_{\Sigma_p})
$$
Here $E_{\Sigma_p}=E|_{\Sigma_p}\otimes\Or_{N^-_X}|_{\Sigma_p}$.
As explained in section~\ref{SS:homotopy} the Morse--Bott--Smale
homotopy, \cf \itemref{BM:iv}, provides a homomorphism of graded complexes
which in view of \itemref{BM:vi} preserves the filtrations:
\begin{equation}\label{E:phi2}
\varphi_\X:C^q_p(X;E)\to C^q(X^+;E)_p
\end{equation}
From \itemref{BM:v} we see that the induced mapping on $E_1$--terms,
\cf Theorem~\ref{T:thom}\itemref{T:thom:ii}, is the isomorphism induced by integration:
%%%%%%%%%%%%%%%%%%%%%%%%%%%%%%%%%%%%%%
%$$
%\begin{CD}
%E_1C^q_p(X;E) @>\varphi_\X>> E_1C^q(X^+;E)_p
%\\
%@| @|
%\\
%H^{q-p}(\Sigma_p;E_{\Sigma_p})
%@>\Int>>
%HC^{q-p}(Y|_{\Sigma_p};E_{\Sigma_p})
%\end{CD}
%$$
%%%%%%%%%%%%%%%%%%%%%%%%%%%%%%%%%%%%%%%%%%%
$$
\xymatrix{
E_1C^q_p(X;E)
\ar[r]^-{\varphi_\X}
\ar@{=}[d]
& 
E_1C^q(X^+;E)_p
\ar@{=}[d]
\\
H^{q-p}(\Sigma_p;E_{\Sigma_p})
\ar[r]^-{\Int}
&
HC^{q-p}(Y|_{\Sigma_p};E_{\Sigma_p})
}
$$

\subsection{Proof of Theorem~\ref{T:main}}\label{SS:proof_main}

The main ingredient for the proof will be the Bismut--Zhang theorem for
Morse--Smale vector fields. We will apply it once for every connected
component $S\subseteq\Sigma$ and once for $M$. Another major input will
be a theorem of Maumary \cite{M69}. This piece of homological algebra 
describes the behavior of torsion with respect to spectral sequences
and generalizes the well known behavior of torsion with respect to short
exact sequences.
The geometric ingredient is the very special Morse--Bott--Smale homotopy
constructed in section~\ref{SS:small_pert}.

We continue to use the notation from section~\ref{SS:small_pert}.
Let $\e\in\Eul_{x_0}(M)$ be the given Euler structure and choose
Euler structures $\e_S\in\Eul_{x_S}(S)$. Choose a one chain $c$ such that
$\e=[-X,c,\{\e_S\}]$ and choose one chains $c_S$ such that
$\e_S=[-Y_S,c_S]$, \cf Lemma~\ref{L:euler}. Set
\begin{equation}\label{E:cc+}
c^+:=c+\sum_{S\subseteq\Sigma}(-)^{\ind(S)}\chi(S)c_S.
\end{equation}
This is an Euler chain for $X^+$ and in view of section~\ref{SS:euler} we know
$\e=[-X^+,c^+]$.

Consider the diagram in figure~\ref{F:CD}. The vertical identifications in
the left column stem from the spectral sequence of the filtered complex
$C^q_p(X;E)$. The corresponding vertical identifications in the
middle column stem from the spectral sequence of the filtered complex
$C^q(X^+;E)_p$. The isomorphisms labeled $c$, $c^+$ and $c_S$
are the ones constructed using parallel transport along these one chains,
see section~\ref{SS:geom_tor}.
We claim that this diagram commutes --- up to sign.

\begin{figure}
\begin{small}
$$
\xymatrix{
&
\det H^*(M;E)
\ar[d]^{\Int_{X^+}}
\ar@/_1pc/[ld]_{\Int_X}
&
(\det E_{x_0})^{\chi(M)}
\ar[l]_{\tau_\an^{E,\e}}
\ar[d]^{c^+}
\ar[dl]^>(0.5){\tau_\geom^{E,\e,X^+}}
\\
\det HC^*(X;E)
\ar[r]^{\varphi_\X}
\ar@{=}[d]
&
\det HC^*(X^+;E)
\ar@{=}[d]
\ar@{=}[r]
&
\det C^*(X^+;E)
\ar@{=}[dd]
\\
\det E_\infty C^*(X;E)
\ar[r]^{\varphi_\X}
\ar@{=}[d]
&
\det E_\infty C^*(X^+;E)
\ar@{=}[d]
\\
\det E_1C^*(X;E)
\ar[r]^{\varphi_\X}
\ar@{=}[d]
&
\det E_1C^*(X^+;E)
\ar@{=}[d]
\ar@{=}[r]
&
\det GC^*(X^+;E)
\ar@{=}[d]
\\
\bigotimes\bigl(\det H^*(S;E_S)\bigr)^\pm
\ar[r]^-{\Int_{Y_S}}
&
\bigotimes\bigl(\det HC^*(Y_S;E_S)\bigr)^\pm
\ar@{=}[r]
&
\bigotimes\bigl(\det C^*(Y_S;E_S)\bigr)^\pm
\\
%(\det E_{x_0})^{\chi(M)}
%\ar[r]_{c}
&
\bigotimes\bigl((\det(E_S)_{x_S})^{\chi(S)}\bigr)^\pm
\ar@/^2pc/[ul]^-{\bigotimes\bigl(\tau_\an^{E_S,\e_S}\bigr)^\pm\quad}
\ar@/_2pc/[ur]_-{\bigotimes(c_S)^\pm}
\ar[u]_-{\bigotimes\bigl(\tau_\geom^{E_S,\e_S}\bigr)^\pm}
\\
&
(\det E_{x_0})^{\chi(M)}
\ar[u]_-{c}
}
$$
\end{small}
\caption{A commutative diagram for the proof of Theorem~\ref{T:main}.
The tensor products are over all connected components $S\subseteq\Sigma$, 
and the powers $\pm$ are abbreviations for $(-)^{\ind(S)}$.}\label{F:CD}
\end{figure}

Indeed, the upper right part commutes because of the Bismut--Zhang theorem applied
to the Morse--Smale vector field $X^+$, but see also Remark~\ref{R:MB_geom_tor}.
The same argument applied to the Morse--Smale vector field $Y_S$ on $S$
(and taking tensor products) shows that the lower half circle commutes. 
The upper left part commutes in
view of section~\ref{SS:homotopy}. The squares involving the mapping
$\varphi_\X$ commute because these mappings are induced from a filtration
preserving map $\varphi_\X:C^q_p(X;E)\to C^q(X^+;E)_p$
and because of the description of $\varphi_\X$ on the $E_1$--term in 
section~\ref{SS:small_pert}.
The big square in the right central part of the diagram commutes in view
of a theorem due to Maumary, see Proposition~\ref{P:spec0} in
appendix~\ref{A:homalg}. The remaining lower right square obviously
commutes.

Note that $\tau_\geom^{E,\e,X}\in\Det_{x_0}(X;E)$
can be seen in the diagram as the composition of
$\otimes\tau_\an^{E_S,\e_S}\circ c$ with the left vertical isomorphisms.
In order to prove Theorem~\ref{T:main} we have to show that this isomorphism
coincides with the composition $\Int_X\circ\tau_\an^{E,\e}$ at the very top.
In view of the commutativity of the diagram it remains
to show that the isomorphism obtained by starting at the very bottom and
going to the upper right corner along the right edge coincides with the
identity on $(\det E_{x_0})^{\chi(M)}$. But this follows immediately from 
\eqref{E:cc+} and the proof is complete.

  \section{Torsion of smooth bundles}
  \label{S:bundles}

\subsection {Bundles} \label{SS:bundles}

Suppose $\pi:M\to N$ is a smooth bundle with $M$ and $N$ closed manifolds,
and suppose $(f',g')$ is a Morse--Smale pair on $N$. Denote by $Y:=-\grad_{g'}(f')$
the corresponding Morse--Smale vector field, and let $\Sigma'$ denote its 
critical manifold which is a discrete set. Define a Morse--Bott function $f$ on $M$ by
$f=f'\circ\pi$. Let us write $\Sigma$ for the critical manifold of $f$.
Clearly $\Sigma=\pi^{-1}(\Sigma')$, and $\ind_\Sigma=\pi^*\ind_{\Sigma'}:\Sigma\to\N_0$. 
Let us also introduce the notation $M_z:=\pi^{-1}(z)$ for the fiber over $z\in N$.

\begin{proposition}\label{P:MBS_on_fb}
\
\begin{enumerate}
\item
There exists a Riemannian metric $g$ on $M$ so that $(f,g)$ 
is a Morse--Bott--Smale pair and so that $X:=-\grad_{g}(f)$ satisfies 
$d_x\pi(X(x))=Y(\pi(x))$ for all $x\in M$.
\item
One can choose $g$ so that $d_x\pi$ restricted to 
$(T_x(M_{\pi(x)}))^\perp$ is an isometry onto $T_{\pi(x)}N$,
and so that the fibers $M_z$ for $z\in\Sigma'$ are isometric. 
\end{enumerate}
\end{proposition}

\begin{proof} 
One chooses a connection in the bundle $\pi$ (\ie for any $x\in M$ a 
subspace $H_x\subset T_xM$ depending smoothly on $x$, so that 
$d_x\pi:H_x\to T_{\pi(x)}N$ is an isomorphism), and a 
smooth family of Riemannian metrics $g_z$, $z\in N$, on the fibers
$M_z$. Clearly $T_xM$ decomposes as $T_xM=H_x\oplus V_x$ with 
$V_x=T_x(M_{\pi(x)})$. Declare $H_x$ and $V_x$ orthogonal and define $g_x$ 
so that $d_x\pi$ is an isometry when restricted $H_x$ and equal to $(g_z)_x$ 
on $V_x$.

To insure that the metric $g$ is compatible with $f$ one chooses a
sufficiently small contractible neighborhood $U\subset N$ of $\Sigma'$ and an identification of 
the bundle $\pi:\pi^{-1}(U)\to U$ with the trivial bundle $\pr_1: U\times F\to U$. 
We chose our connection so that its restriction to 
$\pi^{-1}(U)$ corresponds to the trivial connection in the bundle $\pr_1$, and 
the family of metrics $g_z$, for $z\in U$, corresponds to a constant family of 
Riemannian metrics on $F$. Clearly if $g'$ is compatible with $f'$ then, with the 
above choices, $g$ is compatible with $f$ and satisfies the desired properties.
\end{proof}

Suppose $E$ is a flat vector bundle over $M$.
For $z\in N$ let $E_z:=E|_{M_z}$. Denote by $H^r$ the flat bundle over $N$ whose 
fiber above $z$ is $H^r(M_z;E_z)$.  Note that $H^r$ is equipped with 
an obvious flat connection; any path from $z_1$ to $z_2$ induces an isomorphism from 
$H^r(M_{z_1};E_{z_1})$ to $H^r(M_{z_2};E_{z_2})$ which depends only on the homotopy 
class of the path, and therefore defines a flat connection.
Note that $g_z$ and $\mu_z$ induce, via Hodge theory, a scalar product on
$H^r(M_z;E_z)$ and therefore a fiber metric $\mu^r$ on the bundle $H^r$.

\begin{proposition}\label{P:spec_for_fb}
Suppose $X$ and $Y$ are as in Proposition~\ref{P:MBS_on_fb}. Then:
\begin{enumerate}
\item\label{P:spec_fb_i}
The spectral sequence of the geometric complex $C^*(X;E)$
is isomorphic to the Leray--Serre spectral sequence 
of the bundle $\pi:M\to N$ associated with the cell structure given by $Y$. 
\item\label{P:spec_fb_ii}
$\bigl(E_1C^*(X;E),\delta_1\bigr)$ identifies to 
$\bigoplus_rS^r\bigl(C^*(Y;H^r),\delta\bigr)$.\footnote{Here 
$S^r(C^*,\delta^*)$ denotes the $r$--suspension of the cochain 
complex $(C^*,\delta^*)$, \ie the degree $q$ component of of 
$S^r(C^*,\delta^*)$ is exactly $(C^{q-r},\delta^{q-r})$.}
\end{enumerate}
\end{proposition}

\begin{proof}
Observe that the filtration on $\Omega^*(M;E)$ induced by $X$, see
section~\ref{SS:filt}, induces the Leray--Serre
spectral sequence of the bundle. But the spectral sequence of
$\Omega^*(M;E)$ is isomorphic to the spectral sequence associated to
$C^*(X;E)$, see Remark~\ref{R:specOC} in section~\ref{SS:thom}. 
This proves \itemref{P:spec_fb_i}.

Using $\pi^*(N^-_z)=N^-_{M_z}$ for $z\in\Sigma'$, we get
$\pi^*(\Or_{N^-_z})=\Or_{N^-_{M_z}}$ and hence
$H^*(M_z;E_z\otimes\Or_{N^-_{M_z}})=H^*(M_z;E_z)\otimes\Or_{N^-_z}$.
Now \itemref{P:spec_fb_ii} follows easily from the description of the
differential $\delta_1$ on $E_1C^*(X;E)$ in section~\ref{SS:thom}.
\end{proof}

%%%%%%%%%%%%%%%%%%%%%%%5
%
%1. 
%Denote by $M(r):=\bigcup_{\ind(S)\leq r} W^-_S$ 
%and consider the 
%filtration: 
%$$
%M(0)\subseteq M(1)\subseteq\cdots\subseteq M(r)\subseteq\cdots\subseteq M(n)=M.
%$$ 
%Define $\Omega^*(M;E)_{r}
%:=\{\omega\in\Omega^*(M;E)\mid\omega (x)=0, x\in M(r-1)\}.$ 
%It is not hard to see that, in view of Theorem~\ref{T:comp_W}, $(\Omega^*(M;E)_{r},d^\ast)$ 
%is a cochain subcomplex of $(\Omega^*(M;E), d^*)$  
%whose cohomology is isomorphic to $H(M,M(r-1); E).$    
%
%Define $C^*(X;E)_{r}\equiv C^*(X;E)_{\geq r}:=\bigoplus_{\ind(S)\geq r}
%\Omega^{*-\ind(S)}(S;E_S).$ 
%%which is the same as $C^*(X;E)_{\geq r}$ in Introduction. 
%%and observe that 
%%$C^*(X;E)(r)=C^*(X;E)/C^*(X;E)_{\geq r}$. 
%Note that $Int^*$ induces   
%$Int^*_{r}:\Omega^*(M;E)_r\to C^*(X;E)_{r}$ and these linear maps intertwine $d^\ast$ 
%with $\delta^\ast.$
%%and induces   
%%$Int^*(r):\Omega^*(M;E)(r)\to C^*(X;E)(r)$. 
%By the same arguments as in 
%the proof of Theorem~\ref{T:thom} on concludes that $Int^*_{r}$  
%induces an isomorphism in cohomology. The statement follows immediately from the observation that 
%the filtration $\Omega^*(M;E)_{r}$ induces the Leray--Serre spectral 
%sequence.
%
%2. Part (2) is a actually a tautology.
%%%%%%%%%%%%%%%%%%%%%%%%%%%%%%%%%%%%%%%%%%%

Suppose that $X$ \resp $Y$ are as in Proposition~\ref{P:MBS_on_fb}. 
The isomorphisms between the deRham cohomology and the cohomology
of the geometric complex provided by `integration theory' induces the 
isomorphism
$$
\det\Int^X_E:\det H^*(M;E)\to\det H^*(X;E)
$$
and for any $r\geq0$ the isomorphisms
$$
\det\Int^Y_{H^r}:\det H^*(N;H^r)\to\det H^*(Y;H^r).
$$ 
The last isomorphisms induce the isomorphism 
$$
\det\Int^Y_E:\bigotimes_r\bigl(\det H^*(N;H^r)\bigr)^{\epsilon(r)}\to
\bigotimes_r\bigl(\det H^*(Y;H^r)\bigr)^{\epsilon (r)}
$$
defined by $\bigotimes_r(\det\Int^Y_{H^r})^{\epsilon (r)}$. Here $V^{\epsilon (r)}$  
\resp $\alpha^{\epsilon(r)}$ denote $V$ \resp $\alpha$ if $r$ even and the dual 
$V^*$ \resp $(\alpha^*)^{-1}$ if $r$ is odd.

%The identification of $(E_1C^*(X;E),\delta_1)$ in the spectral sequence of the geometric 
%complex of $X$ to $\bigoplus_{k\geq0}S^k(C^*(Y; H^k),\delta^*)$
%cf Proposition \ref {P:spec_for_fb}  (2.) induces a canonical isomorphism 
%$$\det (\pi, X, Y) :\det H(X;E) \to \otimes_k (\det H(Y; H^k)^{\epsilon (k)}.$$ 

The identification in Proposition~\ref{P:spec_for_fb}\itemref{P:spec_fb_ii}
induces an isomorphism 
$$
E_2C^*(X;E)\to\bigoplus_rH^{*-r}(Y;H^r)
$$ 
and then an isomorphism:
$$
\det(\pi,X,Y):\det H^*(X;E)\to\bigotimes_r\bigl(\det H^*(Y;H^r)\bigr)^{\epsilon (r)}
$$

\begin{proposition}\label{P:det}
The bundle $\pi$ induces a canonical isomorphism
$\det(\pi):\det H^*(M;E)\to\bigotimes_r\bigl(\det H^*(B;H^r)\bigr)^{\epsilon (r)}$ 
so that for $X$ and $Y$ as in Proposition~\ref{P:MBS_on_fb} the following 
diagram is commutative:
$$
\xymatrix{
\det H^*(M;E)
\ar[rr]^-{\det(\pi)}
\ar[d]_{\det\Int^E_X}
&&
\bigotimes_r\bigl(\det H^*(N;H^r)\bigr)^{\epsilon (r)}
\ar[d]^{\det\Int^Y_E}
\\
\det H^*(X;E) 
\ar[rr]^-{\det(\pi,X,Y)}
&&
\bigotimes_r\bigl(\det H^*(Y;H^r)\bigr)^{\epsilon (r)}
}
$$
\end{proposition}

\begin{proof}
Choose a pair $(X,Y)$ as in Proposition~\ref{P:MBS_on_fb}
and  define 
$$
\det(\pi):=(\det\Int^E_Y)^{-1}\circ\det(\pi,X,Y)\circ\Int^E_X.
$$ 
We claim that this does not depend on the choice of $(X,Y)$.
Since given two such pairs $(X_-,Y_-)$ and $(X_+,Y_+)$ one can, by
Proposition~\ref{P:MBS_on_fb}, construct a
Morse--Bott--Smale homotopy which satisfies 
$d_x\pi(\X(x,t))=\mathbb Y(\pi(x),t)$
the arguments in section~\ref{SS:homotopy} imply the independence of the choice of
$(X,Y)$.
\end{proof}

Suppose now that $g$ \resp $g'$ are Riemannian metrics on $M$ \resp $N$.
Recall that $g$ and $\mu$ induce a scalar product on $\det H^*(M;E)$.
Similarly $g'$ and $\mu^r$ induce scalar products on $\det H^*(N;H^r)$. 
Denote by
$$
\log\Vol_\an(\pi):=\log\Vol\det(\pi)
$$ 
and 
$$
\log\Vol_\comb(\pi):=\log\Vol\det(\pi,X,Y).
$$
Proposition~\ref{P:det} implies 
\begin{multline}\label{E:100}
\log\Vol_\an(\pi)-\log\Vol_\comb(\pi)
=\log T_\met^M(E,X,g,\mu)
\\
-\sum_r(-)^r\log T_\met^N(H^r,Y,g',\mu^r)
\end{multline}

\subsection{Torsion for bundles} \label{SS:tor_fb}

As in section~\ref{SS:bundles}, let  
$\pi:M\to N$ be a smooth bundle with $M$ and $N$  closed 
smooth manifolds, $E$ a flat vector bundle over $M$ equipped with a fiber metric
$\mu$, 
and let $g$ \resp $g'$ be a Riemannian metric on $M$ \resp $N$. 
In the previous section we have introduced flat bundles $H^r$
with fiber metrics $\mu^r$.
Denote by $\theta$ \resp $\theta^r$ the closed one forms associated to the fiber 
metric $\mu$ \resp $\mu^r$. Let $X$ and $Y$ be vector fields as in
Proposition~\ref{P:MBS_on_fb}.

Consider the spectral sequence associated with the geometric complex 
$C(X;E)$. For any $k\geq 1$ denote 
$$
\rho^{E,X,g,\mu}_k
:=\frac12\sum_q(-)^{q+1}q\log\det{}'\Delta^q_k 
$$ 
where $\Delta^q_k$ denote the Laplacians in the cochain complexes of 
Euclidean vector spaces $(E_kC^*(X;E),\delta_k)$ and $\det{}'$ denotes the product 
of non-zero eigen values. Note that for $k\geq 1$ $E_kC^*(X;E)$ is a 
finite dimensional vector space, and for $k$ large $\log\det{}'\Delta^q_k=0$. 
Then the calculation of $\rho_k$ for $k\geq1$ is, in principle, a matter of 
finite dimensional linear algebra.

Theorem~\ref{T:intro} applied to $X$ combined with Proposition~\ref{P:spec_for_fb} gives
\begin{multline}\label{E:112}
\log T_\an^M(E,g,\mu)
=\sum_{z\in\Sigma'}(-)^{\ind(z)}\log T_\an^{M_z}(E_z,g_z,\mu_z)
\\
+\sum_r(-)^r\log T_\comb^N(H^r,Y,\mu^r)
+\sum_{k\geq2}\rho^{E,X,g,\mu}_k
\\
+\log T^M_\met(E,X,g,\mu)
+\cR(\theta,-X,g)
\end{multline}
Since $\Sigma'$ is zero dimensional, the definition of $T_\comb^N(\cdot,Y,\cdot)$
does not involve any Riemannian metric on $N$.  
For the metric $g'$ on $N$ the Bismut--Zhang theorem gives 
\begin{multline}\label{E:113}
\log T^N_\comb(H^r,Y,\mu^r)
=\log T^N_\an(H^r,g',\mu^r)
\\
-\log T^N_\met(H^r,Y,g',\mu^r) 
-\cR(\theta^r,-Y,g').
\end{multline}
Combining \eqref{E:100}, \eqref{E:112} and \eqref{E:113} one obtains the following statement

\begin{theorem}\label{T:LST}
Let $\pi:M\to N$ be a smooth fiber bundle with $M$ and $N$ closed manifolds.
Suppose $g$ \resp $g'$ is a Riemannian metric on $M$ \resp $N$,
$X$ and $Y$ are vector fields as in Proposition~\ref {P:MBS_on_fb},
and suppose $E$ is a flat vector bundle over $M$ equipped with a fiber metric $\mu$. 
Then:
\begin{multline*}
\log T^M_\an(E,g,\mu)=
\sum_r(-)^r\log T_\an^N(H^r,g',\mu^r)
\\
+\cR(\theta,-X,g) 
-\sum_r(-)^r\cR(\theta^r,-Y,g') 
+\log\Vol_\an(\pi)
\\
-\log\Vol_\comb(\pi)
+\sum_{k\geq 2}\rho^{E,X,g,\mu}_k
\\
+\sum_{z\in\Sigma'}(-)^{\ind(z)}\log T_\an^{M_z}(E_z,g_z,\mu_z)
\end{multline*}
\end{theorem}

Combining with the metric anomaly for torsion,\footnote{which describes how the 
analytic torsion $T_\an$ changes when one changes the Riemannian metric $g$}
\cf \cite{BZ92}, Theorem~\ref{T:LST} implies the main result of
\cite[Theorems 0.2 and 0.4]{LST98}. 
As in \cite{LST98}, two cases are of particular interest.

\subsubsection*{Bundles with base a sphere.}

Suppose that $N$ is homotopy equivalent to $S^n$. Then for $k\geq1$ 
all $d_k^{p,r}: E^{p,r}_k\to E^{p+k,r-k+1}_k$ are zero but 
$d_n^{0,r}:E_n^{0,r}\to E_n^{n,r-n+1}$. Here we use the notation
$E^{p,r}_k:=E_kC^{p+r}_p(X;E)$ and $d^{p,r}_k$ for the restriction 
$\delta_k:E^{p,r}_k\to E_k^{p+k,r-k+1}$ in order to stay with the 
traditional notation.
Then one obtains the following long exact sequence 
\begin{small}
\begin{equation}\label{E:114}
\cdots 
\to
E^{0,r-1}_n
\xrightarrow{d^{0,r-1}_n}
E^{n,r-n}_n
\xrightarrow{i^{r-n}}
H^r(X;E)
\xrightarrow{p^r}
E^{0,r}_n
\xrightarrow{d^{0,r}_n}
E^{n,r-n+1}_n
\to
\cdots
\end{equation}
\end{small}%
When $H^r(X;E)$ is identified to $H^r(M;E)$ via the integration isomorphism
we refer to the sequence 
\begin{small}
\begin{equation}\label{E:115}
\cdots 
\to
E^{0,r-1}_n
\xrightarrow{d^{0,r-1}_n}
E^{n,r-n}_n
\xrightarrow{i^{r-n}}
H^r(M;E)
\xrightarrow{p^r}
E^{0,r}_n
\xrightarrow{d^{0,r}_n}
E^{n,r-n+1}_n
\to
\cdots
\end{equation}
\end{small}%
as the \emph{Wang exact sequence}.

The Wang exact sequence can be regarded as an acyclic cochain 
complex of finite dimensional Euclidean spaces where each term is equipped 
with the scalar product described above. In particular $H(M;E)$ is equipped
with the scalar product induced from $\mu$ and $g$. 
Denote by $T_W$ the torsion of the complex \eqref{E:115}. Recall that for an 
acyclic cochain complex of Euclidean vector spaces $(C^*,d^*)$ the 
logarithm of the torsion is defined by
$1/2\sum_q(-)^{q+1}q\log\det\Delta^q=\sum_q(-)^q\log\Vol(\overline {d^q})$ 
where $\Delta^q$ is the Laplacian acting on $C^q$. Recall that 
given a linear map $\alpha:V\to W$ between two finite dimensional Euclidean spaces 
one denotes by $\log\Vol(\overline\alpha)$ the quantity, \cf
appendix~\ref{AA:fd_scalar},
$$
\log\Vol(\overline\alpha):=\log\Vol\bigl({\alpha:(\ker\alpha)^\perp\to\alpha(V)}\bigr).
$$

\begin{proposition} \label{P:WES}
One has:
$$
\log T_W= \log T^M_\met(E,X,g,\mu)+\rho^{E,X,g,\mu}_n
$$
\end{proposition}

\begin{proof}
When $H^*(X;E)$ is equipped with a scalar product compatible to the
scalar product on $E^*_\infty$ the torsion of the exact sequence~\eqref{E:114} is
$$
\rho^{E,X,\mu,g}_n=\sum_r(-)^r\log\Vol(\overline{d^{0,r}_n})
$$
simply because $\log\Vol(\overline{i^{n-r}})=0=\log\Vol(\overline{p^r})$.

To see this observe that in the commutative diagram of Euclidean vector
spaces below the left side arrow is surjective, the right side arrow
injective and the lower horizontal row is short exact.
Observe that by construction the scalar product on $E_\infty$ is
induced from the scalar product on $E_n$ and the scalar product on 
$H(X;E)$ is chosen so that it induces the scalar products on 
$E_\infty$ in the lower horizontal row sequence.
$$
\xymatrix{
E^{n,r-n}_n 
\ar[r]^-{i^{r-n}}
\ar[d]
&
H^r(X;E)
\ar[r]^-{p^r}
\ar@{=}[d]
&
E^{n,r-n}_n
\\
E^{n,r-n}_\infty 
\ar[r]^-{i^{n,r-n}_\infty}
&
H^r(X;E)
\ar[r]^{p^{0,r}_\infty}
&
E^{0,r}_\infty.
\ar[u]
}
$$
Give all these we conclude that
\begin{equation*}
\begin{aligned}
\log\Vol(\overline{i^{r-n}}) &= \log\Vol (\overline{i^{n,r-n}_\infty})
\\
\log\Vol(\overline{p^r}) &= \log\Vol(\overline{p^{0,r}_\infty})
\\
\log\Vol(\overline{i^{n,r-n}_\infty})
=&0=\log\Vol(\overline{p^{0,r}_\infty })
\end{aligned}
\end{equation*}
which imply the statement.

Since in the Wang exact sequence $H^r(X;E)$ is replaced by
$H^r(M;E)$ with the scalar product induced by Hodge theory the torsion
$T_W$ of the Wang sequence is corrected by $T_\met^M(E,X,g,\mu)$.
\end{proof}

Consequently, for $n\geq2$, one obtains from \eqref{E:112}: 
\begin{multline*}
\log T_\an^M(E,\mu,g)
=\sum_{z\in\Sigma'}(-)^{\ind(z)}\log T_\an^{M_z}(E_z,g_z,\mu_z)
\\
+\sum_r(-)^r\log T_\comb^{S^n}(H^r,Y,\mu^r)
+\log T_W 
+\cR(\theta,-X,g).
\end{multline*}
For $n=1$ one obtains from Theorem~\ref{T:intro}:
\begin{multline*}
\log T_\an^M(E,\mu,g)
=\sum_{z\in\Sigma'}(-)^{\ind(z)}\log T_\an^{M_z}(E_z,g_z,\mu_z)
\\
+\log T_W 
+\cR(\theta,-X,g).
\end{multline*}

\subsubsection*{Sphere bundles.}

Suppose that the fiber $M_z$ is homotopy equivalent to $S^n$. 
Then  for $k\geq1$ all $d_k^{p,r}:E^{p,r}_k\to E^{p+k,r-k+1}_k$ are 
zero but $d_{n+1}^{p,n}: E_{n+1}^{p,n}\to E_{n+1}^{p+n+1,0}$. From this 
we derive  the following long exact sequence 
\begin{small}
\begin{equation*}%\label{E:122}
\cdots
\to 
E^{p-n-1,n}_{n+1}
\xrightarrow{d^{p-n-1,n}_{n+1}}
E^{p,0}_{n+1}
\xrightarrow{i^p}
H^p(X;E)
\xrightarrow{p^p}
E^{p-n,n}_{n+1}
\xrightarrow{d^{p-n,n}_{n+1}}
E^{p+1,0}_{n+1}
\to
\cdots
\end{equation*}
\end{small}%
When $H^q(X;E)$ is identified to $H^q(M;E)$ via `integration theory' the sequence 
\begin{small}
\begin{equation}\label{E:123}
\cdots
\to 
E^{p-n-1,n}_{n+1}
\xrightarrow{d^{p-n-1,n}_{n+1}}
E^{p,0}_{n+1}
\xrightarrow{i^p}
H^p(M;E)
\xrightarrow{p^p}
E^{p-n,n}_{n+1}
\xrightarrow{d^{p-n,n}_{n+1}}
E^{p+1,0}_{n+1}
\to
\cdots
\end{equation}
\end{small}%
will be referred to as the \emph{Gysin long exact sequence}.

The Gysin long exact sequence can be regarded as an acyclic cochain 
complex of finite dimensional Euclidean spaces where each term is regarded 
with the scalar product 
described above. Denote by $T_G$  the torsion of the complex \eqref{E:123}.

\begin{proposition}\label{P:GES}
One has:
\begin{equation*}
\log T_G=\log T^M_\met(E,X,g,\mu)+\rho^{E,X,g,\mu}_{n+1}
\end{equation*}
\end{proposition}

The proof is identical to the proof of Proposition~\ref{P:WES}.
Consequently we obtain from \eqref{E:112}: 
\begin{multline*}
\log T_\an^M(E,\mu,g)
=\sum_{z\in\Sigma'}(-)^{\ind(z)}\log T_\an^{M_z}(E_z,g_z,\mu_z)
\\
+\log T_\comb^N(H^0,Y,\mu^0)
+(-)^n\log T_\comb^N(H^n,Y,\mu^n)
\\
+\log T_G
+\cR(\theta,-X,g)
\end{multline*}

Note that if $E=\pi^*E'$ is the pull back of a flat bundle $E'$ over $N$
then, as flat bundles over $N$, we have $H^0=H^n=E'$.
Also note that for $n\geq2$ every flat bundle is $E$ over $M$ is the pull back
of a flat bundle $E'$ over $N$.

\begin{appendix}
  \section{Proofs of Theorems~\ref{T:comp_T} and \ref{T:comp_W}}
  \label{A:proof_comp}

The following is a sometimes a simple way to recognize a smooth manifold 
with corners.

\begin{lemma}\label{L:transverse}
If $\mathcal P$ is a smooth manifold with corners, $\mathcal O$ and 
$\mathcal S$ smooth manifolds, $p:\mathcal P\to\mathcal O$ and 
$s:\mathcal S\to\mathcal O$ smooth maps so that $p$ and $s$ are 
transversal then $p^{-1}(s(\mathcal S))$ is a smooth submanifold 
with corners of $\mathcal P$.
Here $p$ is called transversal to $s$ if its restriction 
to each $k$--boundary $\partial_k\mathcal P$ is transversal to $s$.
\end{lemma}

\subsection{Some notations}\label{SS:notations}

\begin{figure}
\includegraphics[scale=0.7]{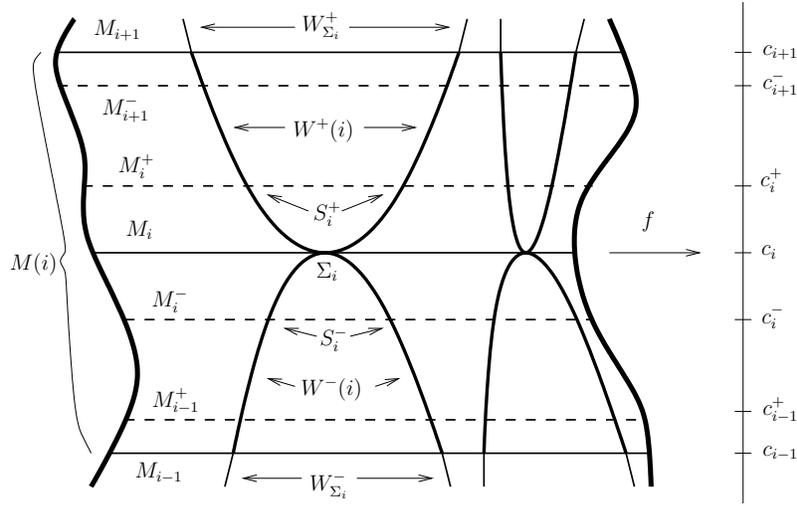}
\caption{Notations.}\label{fig1}
\end{figure}

Let $\cdots<c_i<c_{i+1}<\cdots$ denote the set of all critical values of 
$f$. Choose $\epsilon>0$ sufficiently small, so that
$c_i^+<c_{i+1}^-$, for all $i\in\mathbb Z$, where $c_i^\pm:=c_i\pm\epsilon$.
Define:
\begin{eqnarray*}
M_i      &:=& f^{-1}(c_i)
\\
M_i^\pm  &:=& f^{-1}(c_i^\pm)
\\
M(i)     &:=& f^{-1}(c_{i-1},c_{i+1})
\\
\Sigma_i &:=& \Sigma\cap M_i
\\
W^\pm_{\Sigma_i} &:=& p_\pm^{-1}(\Sigma_i)=\bigsqcup_{S\subseteq\Sigma_i}W^\pm_S
\\
S^\pm_i &:=& W^\pm_{\Sigma_i}\cap M^\pm_i
\\
\mathbf S_i &:=& S^+_i\times_{\Sigma_i}S^-_i
\\
W^\pm(i) &:=& W^\pm_{\Sigma_i}\cap M(i)
\\
\mathbf SW(i) &:=& S^+_i\times_{\Sigma_i}W^-(i)
\end{eqnarray*}
These definitions are to some extent illustrated in Figure~\ref{fig1}.
We have:
\begin{enumerate}
\item
$\mathbf S_i\subseteq M^+_i\times M^-_i$.
\item
$\mathbf SW(i)\subseteq S^+_i\times W^-(i)\subseteq M^+_i\times M(i)$.
\item
$M^\pm_i\subseteq M$ are smooth submanifolds of codimension $1$.
\item
$M(i)$ is an open submanifold of $M$.
\item
$M_i$ is not a manifold, however
$\dot M_i:=M_i\setminus\Sigma_i$ and
$\dot M_i^\pm:=M_i^\pm\setminus S_i^\pm$
are smooth codimension $1$ submanifolds of $M$.
\end{enumerate}

\begin{figure}
\includegraphics[scale=0.7]{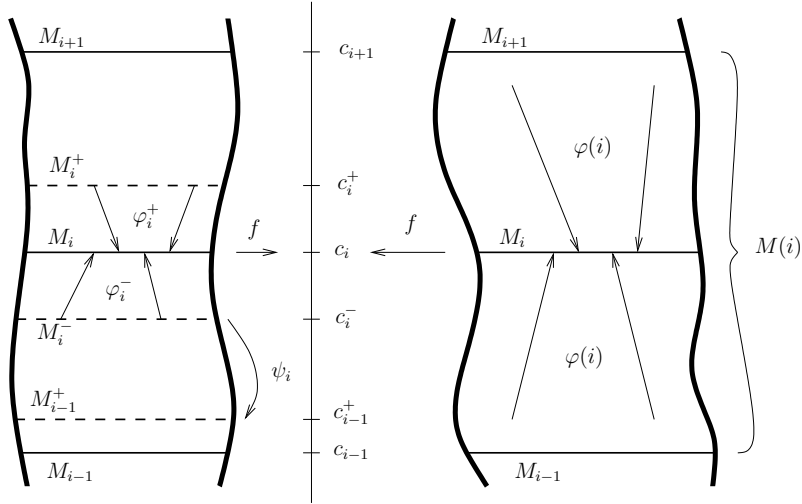}
\caption{Retractions.}\label{fig2}
\end{figure}

Let $\Phi_t$ be the flow associated to the vector field $\frac X{|X|^2}$
on $M\setminus\Sigma$ and consider the diffeomorphisms, see
Figure~\ref{fig2},
$$
\psi_i:M^-_i\to M^+_{i-1},\quad
\psi_i(x):=\Phi_{c_i-\epsilon-(c_{i-1}-\epsilon)}(x)
$$
and
$$
\varphi_i^\pm:\dot M^\pm_i\to\dot M_i,\quad
\varphi_i^\pm(x):=\Phi_{\pm\epsilon}(x),
$$
as well as the submersion
$$
\varphi(i):M(i)\setminus\bigl(W^+(i)\cup W^-(i)\bigr)\to\dot M_i,\quad
\varphi(i)(x):=\Phi_{h(x)-c_i}(x).
$$
The maps $\varphi_i^\pm$ and $\varphi(i)$ extend to
continuous maps
$$
\varphi_i^\pm:M^\pm_i\to M_i
\quad\text{and}\quad
\varphi(i):M(i)\to M_i.
$$
This can be seen from explicit formulae for $\varphi_i^\pm$ and 
$\varphi(i)$ in the neighborhood of $S^\pm_i$ \resp 
$\bigl(W^+(i)\cup W^-(i)\bigr)$, \cf~\cite{BFK}. Define
$$
P_i:=\bigl\{(x,y)\in M_i^+\times M_i^-\bigm|
\varphi^+_i(x)=\varphi^-_i(y)\bigr\},
$$
and denote by $p^\pm_i:P_i\to M_i^\pm$ the canonical projections.
One can verify the following

\begin{lemma} 
$P_i\subseteq M^+_i\times M^-_i$ is a smooth submanifold 
with boundary $\partial P_i$ diffeomorphic to 
$\mathbf S_i\subset M^+_i\times M^-_i$.
Precisely we have:
\begin{itemize}
\item[(P1)]
$p^\pm_i:P_i\setminus\partial P_i\to\dot M_i^\pm$ are diffeomorphisms.
\item[(P2)]
The restriction of $p^+_i\times p^-_i$ to $\partial P_i$
is a diffeomorphism onto $\mathbf S_i$, each $p^\pm_i$ restricted
to $\partial P_i$ identifies with the projection onto $S^\pm_i$.
\end{itemize}
Geometrically $P_i$ represents the space of all trajectories, possibly 
broken, from $M^+_i$ to $M^-_i$.
\end{lemma}

Next we define
$$
Q(i):=\bigl\{(x,y)\in M_i^+\times M(i)\bigm|
\varphi^+_i(x)=\varphi(i)(y)\bigr\}.
$$
Equivalently, $Q(i)$ consists of pairs of points $(x,y)$, $x\in M_i^+$,
$y\in M(i)$, which lie on the same possibly broken trajectory. Moreover
let $l_i:Q(i)\to M^+_i$ and $r_i:Q(i)\to M(i)$ denote the canonical
projections. One can verify the following

\begin{lemma}
$Q(i)\subseteq M^+_i\times M(i)$ is a smooth submanifold with boundary 
$\partial Q(i)$ diffeomorphic to $\mathbf SW(i)\subset M^+_i\times M(i)$. 
Precisely we have:
\begin{itemize}
\item[(Q1)]
$l_i:Q(i)\setminus\partial Q(i)\to\dot M^+_i$ is a smooth bundle with
fiber an open segment.
$r_i:Q(i)\setminus\partial Q(i)\to M(i)\setminus W^-(i)$ is a diffeomorphism.
\item[(Q2)]
The restriction of $l_i\times r_i$ to $\partial Q(i)$
is a diffeomorphism onto $\mathbf SW(i)$, \ie $l_i$ \resp $r_i$ 
restricted to $\partial Q(i)$ identifies with the projection onto
$S^+_i$ \resp $W^-(i)$.
\end{itemize}
\end{lemma}

Since $P_i$ and $Q(i)$ are smooth manifolds with boundaries
\begin{equation}\label{E:Prrk}
\mathcal P_{r,r-k}:=P_r\times P_{r-1}\times\cdots\times P_{r-k}
\end{equation}
and
\begin{equation}\label{E:Qrrk}
\mathcal P_r(r-k):=P_r\times\cdots\times P_{r-k+1}\times Q(r-k)
\end{equation}
are smooth manifolds with corners.

\subsection{Proof of Theorem~\ref{T:comp_T}}

We will equip 
$$
\hat\T(\Sigma_{r+1},\Sigma_{r-k-1})
:=\hat\pi_-^{-1}(\Sigma_{r+1})\cap\hat\pi_+^{-1}(\Sigma_{r-k-1})
=\bigsqcup_{
  {S\subseteq\Sigma_{r+1}}\atop
  {S'\subseteq\Sigma_{r-k-1}}
}
\hat\T(S,S')
$$
with the structure of a smooth manifold with corners. If $k\leq-2$ the 
statement is empty. If $k=-1$ there is nothing to check. So we suppose 
$k\geq 0$. We are going to apply Lemma~\ref{L:transverse}.

Consider the smooth manifold with corners $\mathcal P=\mathcal P_{r,r-k}$
from \eqref{E:Prrk} above, the smooth manifold
$\mathcal O:=\prod^{r-k}_{i=r}(M_i^+\times M_i^-)$ and the smooth manifold
$\mathcal S:=S^-_{r+1}\times M^-_r\times\cdots
\times M^-_{r-k+1}\times S^+_{r-k-1}$.
In order to define the maps $p$ and $s$ we consider
$$
\omega_i:M^-_i\to M^-_i\times M^+_{i+1},\quad
\omega_i(x):=(x,\psi_i(x))
$$
and
$$
\tilde p_i:P_i\to M^+_i\times M^-_i,\quad
\tilde p_i(y):=\bigl(p^+_i(y),p^-_i(y)\bigr).
$$
We also denote by $\alpha:S^-_{r+1}\to M^+_r$ \resp 
$\beta:S^+_{r-k-1}\to M^-_{r-k}$ the restriction of $\psi_{r+1}$ \resp 
$\psi_{r-k}^{-1}$ to $S_{r+1}^-$ \resp $S_{r-k-1}^+$.
Finally we set, 
$$
s:=\alpha\times\omega_r\times\cdots\times\omega _{r-k+1}\times\beta:
\mathcal S\to\mathcal O
$$
and
$$
p:=\tilde p_r\times\cdots\times\tilde p_{r-k}:
\mathcal P\to\mathcal O.
$$

The verification of the transversally of $p$ and $s$ follows easily from
(P1), (P2) and the Morse--Smale condition and is done exactly as in 
\cite{BH01}. It is easy to see that $p^{-1}(s(\mathcal S))$ identifies to
$\hat{\mathcal T}(\Sigma_{r+1},\Sigma_{r-k-1})$ as topological spaces and
we leave this verification to the reader. The rest is straight forward.

\subsection{Proof of Theorem~\ref{T:comp_W}}

We are going equip the compact topological space 
$X:=\hat W^-_{\Sigma_{r+1}}$ with the structure of
a smooth manifold with corners. Let us write
$\hat i:=\hat i^-_{\Sigma_{r+1}}:X\to M$. For any integer $k$ 
denote by $X(r-k):=\hat i^{-1}(M(r-k))$. First we will 
put the structure of smooth manifold with 
corners on $X(r-k)$ so that the restriction of $\hat i$
is smooth. Second we check that $X(r-k)$ and 
$X(r-k')$ induce on the intersection $X(r-k)\cap X(r-k')$ the same 
smooth structure. These facts imply that $X$ has a canonical structure of 
smooth manifold with corners and that $\hat i$ is a smooth map.

To accomplish first step we proceed in exactly the same way as in the
proof of Theorem~\ref{T:comp_T}. Consider the smooth manifold with 
corners $\mathcal P:=\mathcal P_r(r-k)$ from \eqref{E:Qrrk}, the smooth 
manifold $\mathcal O:=\prod_{i=r}^{r-k+1}(M_i^+\times M_i^-)\times M^+_{r-k}$
and the smooth manifold 
$\mathcal S:=S^-_{r+1}\times M^-_r\times\cdots\times M^-_{r-k+1}$.
Define
$$
p:=\tilde p_r\times\cdots\times\tilde p_{r-k+1}\times l_{r-k}:
\mathcal P\to\mathcal O
$$
and
$$
s:=\alpha\times\omega_r\times\cdots\times\omega_{r-k+1}:\mathcal S\to\mathcal O.
$$
Again we will apply Lemma~\ref{L:transverse}.
The verification of the transversality follows from (P1), (P2), (Q1), (Q2)
and the Morse--Smale condition, exactly as explained in \cite{BH01}. 
It is easy to see that $p^{-1}(s(\mathcal S))$
identifies to $X(r-k)$. The rest is straight forward.

  \section{Some homological algebra}
  \label{A:homalg}

\subsection{Finite dimensional complexes}\label{AA:det}

All vector spaces in this section are supposed to be finite dimensional
and over a field $\K$ of characteristics zero.
For such a vector space $V$ we define its
\emph{determinant} by $\det V:=\Lambda^{\dim(V)}V$.
If we have a short exact sequence of vector spaces
$0\to U\to V\to W\to0$ we obtain a canonic isomorphism
$\det V=\det U\otimes\det W$. This is constructed using a section 
$W\to V$ to get an isomorphism $V=U\oplus W$. This induces an isomorphism
$\det V=\det U\otimes\det W$ which does not depend on the section.

If $V=V^\even\oplus V^\odd$ is a $\Z_2$--graded vector space
we define its determinant by $\det V:=\det V^\even\otimes\det(V^\odd)^*$.
Here $(V^\odd)^*$ denotes the dual of $V^\odd$. Note that we have a natural
identification $\det((V^\odd)^*)=(\det V^\odd)^*$.
Again, a short exact sequence $0\to U\to V\to W\to0$
of $\Z_2$--graded vector spaces gives rise to a natural isomorphism
$\det V=\det U\otimes\det W$ by combining the isomorphisms
for the even and the odd part.

Suppose $C=C^\even\oplus C^\odd$ is a $\Z_2$--graded complex
with a degree differential of odd degree $d:C^\even\to C^\odd$ and
$d:C^\odd\to C^\even$.
Then we get three $\Z_2$--graded vector spaces;
$Z:=\ker d$, $B:=\img d$ and $HC:=Z/B$.
Moreover we get two short exact sequences of $\Z_2$--graded vector spaces
$0\to B\to Z\to HC\to0$ and
$0\to Z\to C\overset d\to B^\op\to0$. Here $B^\op$ denotes the
$\Z_2$--graded vector space $B$ equipped with the opposite grading,
$(B^\op)^\even:=B^\odd$ and $(B^\op)^\odd:=B^\even$. Note that we have a
natural isomorphism $\det B^\op=(\det B)^*$.
The two short exact sequences provide natural isomorphisms
$$
\det Z=\det B\otimes\det HC
\quad\text{and}\quad
\det C=\det Z\otimes\det B^\op.
$$
Together with $\det B\otimes\det B^\op=\K$ they give rise to a
natural isomorphism 
\begin{equation}\label{E:CdetHC}
\det C=\det HC.
\end{equation}

Suppose we have a short exact sequence of $\Z_2$--graded complexes
$0\to C_0\to C_1\to C_2\to0$. Its long exact sequence in cohomology
$$
\xymatrix{
H^\even C_0 
\ar[r]
&
H^\even C_1 
\ar[r]
& 
H^\even C_2
\ar[d]
\\
H^\odd C_2 
\ar[u]
&
H^\odd C_1 
\ar[l]
&
H^\odd C_0
\ar[l]
}
$$
can be used to make the $\Z_2$--graded vector space
$HC_0\oplus(HC_1)^\op\oplus HC_2$
a $\Z_2$--graded complex with vanishing cohomology. From \eqref{E:CdetHC} we
thus obtain a natural isomorphism:
\begin{eqnarray}\label{E:les}
\det HC_0\otimes\det HC_2=\det HC_1
\end{eqnarray}

The following is a reformulation of a well known fact, which for instance 
can be found in \cite[Theorem~3.2]{Mi66}.

\begin{proposition}\label{P:det_ses}
Suppose $0\to C_0\to C_1\to C_2\to 0$ is a short exact sequence of
$\Z_2$--graded complexes. Then the following diagram of canonic
isomorphisms commutes up to sign:
$$
\xymatrix{
\det C_0\otimes\det C_2 
\ar@{=}[r]
&
\det C_1
\ar@{=}[d]
\\
\det HC_0\otimes\det HC_2 
\ar@{=}[u]
& 
\det HC_1
\ar@{=}[l]
}
$$
\end{proposition}

Suppose $C$ is a $\Z_2$--graded complex equipped with a filtration
$$
\cdots\supseteq C_p\supseteq C_{p+1}\supseteq\cdots
$$
which is preserved by the differential. We assume $C_p=0$ for $p$
sufficiently large and $C_p=C$ for $p$ sufficiently small. Moreover we
assume that the filtration is compatible with the grading, that is 
$C_p=C_p^\even\oplus C_p^\odd$ where $C_p^\even:=C_p\cap C^\even$ and
$C_p^\odd:=C_p\cap C^\odd$. In this situation we have a spectral sequence
$(E_kC,\delta_k)$ converging to $HC$. Every $E_kC$ is a $\Z_2$--graded
complex with differential $\delta_k$. For the $E_1$--term we have 
$E_1C=HGC$ as $\Z_2$--graded vector spaces. Here
$GC:=\bigoplus_pC_p/C_{p+1}$ denotes the associated graded of $C$
considered as $\Z_2$--graded
complex, and $HGC$ its cohomology, a $\Z_2$--graded vector space. We will
not consider the $\Z$--grading $GC$ and $HGC$ inherit from the filtration.
Applying  \eqref{E:CdetHC} we obtain a canonic isomorphism 
$\det C=\det GC=\det HGC=\det E_1C$. The cohomology $HC$
inherits the structure of a $\Z_2$--graded vector space with filtration
from $C$. Let $GHC$ denote the associated graded considered as
$\Z_2$--graded vector space. Again, we will not consider the $\Z$--grading
$GHC$ inherits. Since the spectral sequence converges we
have a canonic isomorphism $GHC=E_\infty C$ which provides us with a canonic
isomorphism $\det C=\det HC=\det GHC=\det E_\infty C$. Since
$HE_kC=E_{k+1}C$ we also have a canonic isomorphism $\det
E_k=\det E_{k+1}C$ for every $k\geq1$. Also note that $E_k=E_\infty$ for
$k$ sufficiently large.

The following can be found in \cite[Theorem 1.3]{F92}, but see also
\cite[Th\'eor\`em 4.4]{M69} for a more general result.

\begin{proposition}\label{P:spec0}
Let $C$ be a filtered $\Z_2$--graded complex. Then the following diagram
of canonic isomorphism commutes up to sign:
$$
\xymatrix{
\det HGC 
\ar@{=}[r]
&
\det GC 
\ar@{=}[r]
&
\det C 
\ar@{=}[r]
&
\det HC 
\ar@{=}[r]
&
\det GHC
\ar@{=}[d]
\\
\det E_1C 
\ar@{=}[u]
&
\det E_2C
\ar@{=}[l]
&
\det E_3C
\ar@{=}[l]
&
\cdots 
\ar@{=}[l]
&
\det E_\infty C
\ar@{=}[l]
}
$$
\end{proposition}

\subsection{Finite dimensional complexes with scalar product}
\label{AA:fd_scalar}

All vector spaces (complexes) in this section are suppose to be 
finite dimensional and over $\R$. We will work with $\Z_2$--graded complexes
and all determinant lines are understood in the graded sense. If a
$\Z_2$--graded complex $C=C^\even\oplus C^\odd$ comes with a
scalar product $g$ we will assume $g=g^\even\oplus g^\odd$, \ie the
even and the odd part are orthogonal. In
such a situation $B=\img d$, $Z=\ker d$ and $HC=Z/B$ are
$\Z_2$--graded vector spaces with scalar products. We let 
$||\cdot||_{\det C}$, $||\cdot||_{\det B}$, $||\cdot||_{\det Z}$ and
$||\cdot||_{\det HC}$ denote the induced norms
on $\det C$, $\det B$, $\det Z$ and $\det HC$, respectively. Moreover let
$d^\sharp$ denote the adjoint of the differential $d$ and $\Delta=d^\sharp
d+dd^\sharp$ the Laplacian.

For a linear isomorphism between finite dimensional Euclidean spaces 
$\alpha:V\to W$ we set $\Vol\alpha:=\sqrt{\det(\alpha^\sharp\alpha)}$. For an
arbitrary linear mapping $\alpha:V\to W$ between Euclidean spaces we
let $\overline\alpha:(\ker\alpha)^\perp\to\img\alpha$ denote its 
restriction and thus have a volume $\Vol(\overline\alpha)$. Moreover,
if $\beta:U\to U$ is a linear mapping between finite dimensional 
vector spaces we write $\det{}'\beta$ for the product of non-zero 
eigen values. Clearly $\Vol(\overline\alpha)=\sqrt{\det{}'(\alpha^\sharp\alpha)}$.
For a $\Z_2$--graded complex $C$ with scalar product define:
$$
T_C:=\frac{\Vol\bigl(\overline{d:C^\even\to01z C^\odd}\bigr)}
{\Vol\bigl(\overline{d:C^\odd\to C^\even}\bigr)}
=\frac{\sqrt{\det{}'(d^\sharp d:C^\even\to C^\even)}}
{\sqrt{\det{}'(d^\sharp d:C^\odd\to C^\odd)}}
$$

\begin{remark}
Suppose we have a $\Z$--graded complex $C=\bigoplus_qC^q$ with
scalar product. It gives rise to a $\Z_2$--graded complex with scalar
product in the usual way. Let $\Delta^q$ denote the Laplacian acting on
$C^q$. A two line computation using
$$
\det{}'(\Delta^q)=\det{}'(d^\sharp d:C^q\to C^q)\det{}'(dd^\sharp:C^q\to C^q)
$$
and
$$
\det{}'(dd^\sharp:C^q\to C^q)=\det{}'(d^\sharp d:C^{q-1}\to C^{q-1})
$$
shows
$$
\log T_C=\frac12\sum_q(-)^{q+1}q\log\det{}'\Delta^q
$$
However, for this to make sense it is essential to have a $\Z$--grading at
hand.
\end{remark}

\begin{proposition}\label{P:lap_tor}
Let $C$ be a $\Z_2$--graded complex with scalar product. Then the canonic
isomorphism $\det C=\det HC$ maps $||\cdot||_C$ to $T_C||\cdot||_{HC}$.
\end{proposition}

\begin{proof}
The short exact sequence $0\to B\to Z\to HC\to0$ of $\Z_2$--graded
vector spaces induces an isomorphism $\det Z=\det B\otimes\det HC$ which
clearly maps $||\cdot||_{\det Z}$ to 
$||\cdot||_{\det B}\otimes||\cdot||_{\det HC}$. The short exact sequence
$0\to Z\to C\to B^\op\to0$ of $\Z_2$--graded vector spaces induces an
isomorphism $\det C=\det Z\otimes\det B^\op$ which maps $||\cdot||_{\det C}$
to $T_C||\cdot||_{\det Z}\otimes||\cdot||_{\det B^\op}$. Putting this
together we find $||\cdot||_{\det C}=T_C||\cdot||_{\det HC}$ via
$\det C=\det HC$.
\end{proof}

\begin{proposition}\label{P:ses}
Suppose $0\to C_0\to C_1\to C_2\to0$ is a short exact sequence of 
$\Z_2$--graded complexes with compatible scalar products, \ie
$C_0\to C_1$ is an isometry onto its image and $C_2$ carries the quotient
scalar product. Then the canonic isomorphism
$$
\det HC_1
=\det HC_0\otimes\det HC_2
$$
maps $T_{C_1}||\cdot||_{HC_1}$ to 
$T_{C_0}||\cdot||_{HC_0}\otimes T_{C_2}||\cdot||_{HC_2}$.
\end{proposition}

\begin{proof}
The compatibility condition for the scalar products guarantees that the 
isomorphism $\det C_0\otimes\det C_2=\det C_1$ maps 
$||\cdot||_{\det C_0}\otimes||\cdot||_{\det C_2}$ to
$||\cdot||_{\det C_1}$. Now the statement follows from
Proposition~\ref{P:det_ses} and Proposition~\ref{P:lap_tor}.
\end{proof}

Let $C$ be a filtered $\Z_2$--graded complex as considered
in section~\ref{AA:det} and suppose $C$ is equipped with a scalar product.
Then $HC$, $GC$ and $HGC$ inherit scalar products in an obvious way.

\begin{proposition}\label{P:spec}
Suppose $C$ is a filtered $\Z_2$--graded complex with scalar product. 
Then the canonic isomorphism
$$
\det HC=\det GHC=\det E_\infty C=\cdots=\det E_1C=\det HGC
$$
maps $T_C||\cdot||_{HC}$ to $T_{GC}||\cdot||_{HGC}$.
\end{proposition}

\begin{proof}
This follows immediately from Proposition~\ref{P:spec0},
Proposition~\ref{P:lap_tor} applied to $C$ and $GC$ as well as the fact that
the natural isomorphism $\det C=\det GC$ maps $||\cdot||_{\det C}$ to
$||\cdot||_{\det GC}$.
\end{proof}

%  \section{Examples and applications}
%  \label{A:ex}
%  \input{ex.tex}
\end{appendix}


\begin{thebibliography}{BFK99}
     


\bibitem[AB95]{AB95}
    D.M.~Austin and P.J.~Braam,
    \emph{Morse--Bott theory and equivariant cohomology}. 
    The Floer memorial volume, 123--183.
    Progr. Math. \textbf{133}, 
    Birkh\"auser, Basel, 1995.


\bibitem[BZ92]{BZ92}
    J.M.~Bismut and W.~Zhang,
    \emph{An extension of a theorem by Cheeger and M\"uller}.
    With an appendix by F.~Laudenbach.
    Astérisque \textbf{205}(1992).


\bibitem[BZ94]{BZ94}
    J.M. Bismut and W. Zhang,
    \emph{Milnor and Ray--Singer metrics on the equivariant determinant 
          of a flat vector bundle}.
    Geom. Funct. Anal. \textbf{4}(1994), 136--212.


\bibitem[BT82]{BT82}
    R.~Bott and L.W.~Tu, 
    \emph{Differential forms in algebraic topology}.
    Graduate Texts in Mathematics \textbf{82}. 
    Springer--Verlag, 
    New York--Berlin, 1982.


\bibitem[B99]{B99}
    D.~Burghelea,
    \emph{Removing metric anomalies from Ray--Singer torsion},
    Lett. Math. Phys. \textbf{47}(1999), 149--158.


\bibitem[BFK99]{BFK99}
  D. Burghelea, L. Friedlander and T. Kappeler,
  \emph{Torsions for manifolds with boundary and glueing formulas}.
  Math. Nachr. \textbf{208}(1999), 31--91.
  

\bibitem[BFK01]{BFK01}
  D. Burghelea, L. Friedlander and T. Kappeler,
  \emph{Relative torsion}.
  Commun. Contemp. Math. \textbf{3}(2001), 15--85.


\bibitem[BFK]{BFK}
  D.~Burghelea, L.~Friedlander and T.~Kappeler,
  \emph{Elementary Morse Theory},
  Lectures Notes in preparation.


\bibitem[BH01]{BH01}
    D.~Burghelea and S.~Haller,  
    \emph{On the topology and analysis of a closed one form.~I
          (Novikov's theory revisited)}.  
    Essays on geometry and related topics, 133--175.
    Monogr. Enseign. Math. \textbf{38}, 
    Enseignement Math., Geneva, 2001.


\bibitem[BH]{BH}
    D.~Burghelea and S.~Haller,  
    \emph{A Riemannian invariant, Euler structures and some topological
          applications},
    preprint {\tt math.DG/0310154}.


\bibitem[BH04]{BH04}
    \dots


\bibitem[C77]{Ch77}
    J.~Cheeger, 
    \emph{Analytic torsion and Reidemeister torsion},
    Proc. Nat. Acad. Sci. U.S.A. \textbf{74}(1977), 2651--2654.


\bibitem[C79]{Ch79}
    J.~Cheeger,
    \emph{Analytic torsion and the heat equation},
    Ann. of Math. \textbf{109}(1979), 259--322.


\bibitem[D95]{D95}
    A. Dold,
    \emph{Lectures on algebraic topology}.
    Reprint of the 1972 edition.
    Classics in Mathematics.
    Springer--Verlag, Berlin, 1995.


\bibitem[F92]{F92}
    D.S.~Freed, 
    \emph{Reidemeister torsion, spectral sequences, and Brieskorn spheres},
    J. Reine Angew. Math. \textbf{429}(1992), 75--89.


\bibitem[LST98]{LST98}
    W.~L\"uck, T.~Schick and T.~Thielmann,
    \emph{Torsion and fibrations},
    J. Reine Angew. Math. \textbf{498}(1998), 1--33.


\bibitem[M69]{M69}
    S.~Maumary, 
    \emph{Contributions \`a la th\'eorie du type simple d'homotopie},
    Comment. Math. Helv. \textbf{44}(1969), 410--437.


\bibitem[M66]{Mi66}
    J.~Milnor,
    \emph{Whitehead torsion}, 
    Bull. Amer. Math. Soc. \textbf{72}(1966), 358--426.


\bibitem[M78]{Mu78}
    W.~M\"uller, 
    \emph{Analytic torsion and $R$-torsion of Riemannian manifolds},
    Adv. in Math. \textbf{28}(1978), 233--305.


\bibitem[RS71]{RS71}
    D.B.~Ray and I.M.~Singer,
    \emph{$R$--torsion and the Laplacian on Riemannian manifolds},
    Adv. in Math. \textbf{7}(1971), 145--210.


\bibitem[T90]{Tu90}
    V.~Turaev,
    \emph{Euler structures, non-singular vector fields and 
          Reidemeister-type torsions},
    Math. USSR--Izv. {\bf 34}(1990), 627--662.


\bibitem[T02]{Tu02}
    V.~Turaev,
    \emph{Torsions of $3$--dimensional manifolds}. 
    Progress in Mathematics \textbf{208}.
    Birkh\"auser Verlag, Basel, 2002.




  \end{thebibliography}
\end{document}